\documentclass[11pt]{amsart}

\usepackage{amsmath,amsthm,amsfonts,amssymb,bm,graphicx}

\usepackage%[backref=page]
{hyperref}
\hypersetup{colorlinks=true,citecolor=blue,linkcolor=blue,urlcolor=blue,
pdfstartview=FitH }

\theoremstyle{plain}
\newtheorem{theorem}{Theorem}
\newtheorem{cor}[theorem]{Corollary}
\newtheorem{lemma}{Lemma}

\theoremstyle{definition}

\renewcommand{\geq}{\geqslant}
\renewcommand{\leq}{\leqslant}

\newcommand{\dd}{{\rm d}}
\numberwithin{equation}{section}

\makeatletter
\def\Ddots{\mathinner{\mkern1mu\raise\p@
\vbox{\kern7\p@\hbox{.}}\mkern2mu
\raise4\p@\hbox{.}\mkern2mu\raise7\p@\hbox{.}\mkern1mu}}
\makeatother
\makeatletter
\DeclareRobustCommand\widecheck[1]{{\mathpalette\@widecheck{#1}}}
\def\@widecheck#1#2{%
    \setbox\z@\hbox{\m@th$#1#2$}%
    \setbox\tw@\hbox{\m@th$#1%
       \widehat{%
          \vrule\@width\z@\@height\ht\z@
          \vrule\@height\z@\@width\wd\z@}$}%
    \dp\tw@-\ht\z@
    \@tempdima\ht\z@ \advance\@tempdima2\ht\tw@ \divide\@tempdima\thr@@
    \setbox\tw@\hbox{%
       \raise\@tempdima\hbox{\scalebox{1}[-1]{\lower\@tempdima\box
\tw@}}}%
    {\ooalign{\box\tw@ \cr \box\z@}}}
\makeatother

\begin{document}

\author{Valentin Blomer}
\address{Mathematisches Institut, Endenicher Allee 60, 53115 Bonn, Germany} \email{blomer@math.uni-bonn.de}

 \title{Density theorems for  ${\rm GL}(n)$}

\thanks{Author partially supported by DFG grant BL 915/2-2}

\keywords{exceptional eigenvalues, Kloosterman sums, Kuznetsov formula, Ramanujan conjecture}

\begin{abstract} Strong bounds -- going beyond Sarnak's density hypothesis -- are obtained for the number of automorphic forms for the group $\Gamma_0(q) \subseteq {\rm SL}(n, \Bbb{Z})$ violating the Ramanujan conjecture at any given unramified   place.   The proof is based on a relative trace formula of Kuznetsov type and  best-possible bounds for certain Kloosterman sums for ${\rm GL}(n)$. Further applications are given. 
 \end{abstract}

\subjclass[2010]{Primary 11F72, 11L05}

\setcounter{tocdepth}{2}  \maketitle 

\maketitle

\section{Introduction}

The concept of families is a very fruitful one in number theory and in particular in the context of automorphic forms. It allows us to study asymptotic properties and  has recently been put on some formal ground  in \cite{SST, Sa2}. On the conceptual side it dampens irregularities of individual members (that may exist or whose non-existence we are unable to prove) and allows statistical concepts and deformation techniques to investigate properties within an ensemble. On the methodological it enables us to use strong analytic tools such as various types of trace formulae.

One of the key conjectures in the field of automorphic forms is the Ramanujan conjecture: cuspidal automorphic representations of the group ${\rm GL}(n)$ over a number field $F$ are tempered (see \cite{BlBr} for a survey). Even for $n = 2$ this appears to be far out of reach, and as a substitute one considers two types of approximations. On the one hand one can measure the worst case scenario, i.e. the largest distance from the tempered spectrum of an individual member in a family. On the other hand one can try to bound the number of members in a family violating the conjecture relative to the amount by which they violate the conjecture. This is a density result which is a familiar concept from the theory Dirichlet $L$-functions: although the Riemann hypothesis is far out of reach, we have good bounds for the number $N(\sigma,T, Q)$ of zeros with real part $\geq \sigma$ and height $\leq T$ of Dirichlet $L$-functions with conductor $q \leq  Q$ (see e.g.\ \cite[Section 10]{IK}). The arithmetic reformulation of this is the Bombieri-Vinogradov theorem which roughly states that primes $\leq x$ are equidistributed in ``almost all'' residue classes modulo $q \leq x^{1/2 + o(1)}$  (similarly, ``almost all'' short intervals contain primes). In many applications this serves as a good substitute for the Riemann hypothesis.

In this note we want to consider the automorphic analogue for the family of automorphic forms for the group $\Gamma_0(q) \subseteq {\rm SL}_n(\Bbb{Z})$ of matrices whose lowest row is congruent to $(0, \ldots, 0, \ast)$ modulo $q$. This is a very natural family as it contains precisely the automorphic forms of conductor dividing $q$ \cite{JPSS}. Let us fix a place $v$ of $\Bbb{Q}$, and for an automorphic form $\pi$ let us denote by $\mu_{\pi}(v) = (\mu_{\pi}(v, 1), \ldots, \mu_{\pi}(v, n))$ its local spectral parameter (each entry viewed modulo  $\frac{2\pi i}{\log p} \Bbb{Z}$  if $v = p$ is a prime). %Here and in the following we let $W$ be the Weyl group of permutation matrices and we identify a permutation matrix $w = (w_{ij}) \in W$ with the  permutation $ i \mapsto j$ for $w_{ij} = 1$. 
Write 
\begin{equation}\label{sigma}
  \sigma_{\pi}(v) = \max_j|\Re \mu_{\pi}(v, j)|.
 \end{equation} 
    The representation $\pi$ is tempered at $v$ if $\sigma_{\pi}(v) = 0$, and the size of $\sigma_{\pi}(v)$ measures how far $\pi$ is from being tempered at $v$.  An example of a non-tempered representation is the trivial representation which satisfies $\sigma_{\text{triv}}(v) = (n-1)/2$ for every $v$. %Fix some bounded, measureable set $\Omega \subseteq \mathfrak{a}_{\Bbb{C}}^{\ast}$ of archimedean spectral parameters and let $\mathcal{F}(q, \Omega)$ denote the family of automorphic 
For a finite  family $\mathcal{F}$ of automorphic representations for ${\rm GL}(n)$ and $\sigma \geq 0$ we define
$$N_{ v}(\sigma, \mathcal{F}) = |\{\pi \in \mathcal{F} \mid  \sigma_{\pi}(v) \geq \sigma\}|.$$
 We have trivially $N_{ v}(0, \mathcal{F}) = |\mathcal{F}|$, and if the trivial representation is contained in $\mathcal{F}$, we have $N_{ v}((n-1)/2, \mathcal{F}) \geq 1$.  One may hope to   be able to interpolate linearly between these two extreme cases:
\begin{equation}\label{density}
N_{ v}(\sigma, \mathcal{F}) \ll_{v, \varepsilon} |\mathcal{F}|^{1 - \frac{\sigma}{a}+ \varepsilon}
\end{equation}
for arbitrarily small $\varepsilon > 0$ with 
\begin{equation}\label{c}
 a = (n-1)/2.
  \end{equation}
This is precisely Sarnak's density hypothesis \cite[p.\ 465]{Sa1} stated there in the context of groups $G$ of real rank 1, the principal congruence subgroup $\Gamma(q) = \{\gamma \in G(\Bbb{Z}) \mid \gamma \equiv \text{id} \, (\text{mod } q)\}$ and $v = \infty$. For families of large level, Sarnak's density hypothesis  has recently attracted  interest in the context of lifting matrices modulo $q$ \cite{Sa3} and the almost diameter of Ramanujan complexes, and for families with growing infinitesimal character in the context of Golden Gates and quantum computing \cite{Sa6, PS}. In each of these cases it is not a spectral gap that is needed, but a certain kind of density result. 

The shape of the bound \eqref{density} -- \eqref{c}   bears a certain similarity to the convexity bound for $L$-functions in the Selberg class in the critical strip. Unlike the convexity bound for $L$-functions, \eqref{density}  -- \eqref{c} is in general a very deep result that is completely open for general groups and families. On the other hand, it is a priori not impossible to even obtain ``subconvexity'', i.e.\ a proof of \eqref{density} with a constant $a < (n-1)/2$, if the trivial representation is not in $\mathcal{F}$. The Arthur-Selberg trace formula is usually not sensitive to whether the trivial representation is counted or not, but the Kuznetsov formula can be a versatile tool if no residual spectrum is involved.

For the group ${\rm GL}(2)$ there exist strong density results for many automorphic families, for instance by Sarnak \cite{Sa}, Iwaniec \cite{Iw}, Huxley \cite{Hu}, Blomer-Buttcane-Raulf \cite{BBR}, also in number field versions \cite{BM1, BM2} and for general real rank 1 groups \cite{SX, HK}. Various results are also available for $ {\rm GL}(3)$, see e.g. \cite{Bl, BBM, BBR}. For higher rank groups the deep analysis of the Arthur-Selberg trace formula of Matz-Templier \cite{MT} and Finis-Matz \cite{FM} provides as by-products some density results   for the family of Maa{\ss} forms   of Laplace eigenvalue up to height $T$  and fixed level. The value of $a$ is however much larger than \eqref{c} for $n > 2$ (at least quadratic in $n$), so that even the ``convexity bound'' cannot be obtained. 

In the present paper we consider the family $\mathcal{F}_I(q)$ of cuspidal automorphic representations generated by   Maa{\ss}  forms for the group $\Gamma_0(q) \subseteq {\rm SL}_n(\Bbb{Z})$  for a large prime $q$  and Laplace eigenvalue $\lambda$ in a fixed   interval $I$. If $I$ is not too small, we have   $|\mathcal{F}_I(q)| \asymp_I  q^{n-1}$. For this family and any place $v \not= q$ of $\Bbb{Q}$, we %prove Sarnak's density conjecture \eqref{density}. If in addition $q$ is prime,  we  
go \emph{beyond} the density hypothesis and obtain strong ``subconvexity''   with a value of $$a = (n-1)/4,$$
 which is halfway between \eqref{c} and the Ramanujan conjecture. %Note that $|\mathcal{F}_I(q)| = q^{n-1 + o_I(1)}$. 

 \begin{theorem}\label{thm1} Let $n \geq 3$, $q$ a prime, $v$ be a place of $\Bbb{Q}$ different from $q$,  $I \subseteq [0, \infty)$ a fixed interval, $\varepsilon > 0$, and $\sigma \geq  0$. Then
 %$$N_v(\sigma, \mathcal{F}(q, I)) \ll_{I, v, \varepsilon}  q^{n-1 - 2\sigma + \varepsilon}$$
 %If $q$ is prime, then  
 $$N_v(\sigma, \mathcal{F}_I(q)) \ll_{I, v, n, \varepsilon}  q^{n-1 - 4\sigma + \varepsilon}.$$
 \end{theorem}
 
 Of course, by \cite{LRS} we know that $N_v(\sigma, \mathcal{F}_I(q)) = 0$ for $\sigma \geq 1/2 - 1/(n^2 + 1)$, but for $0 < \sigma < 1/2 - 1/(n^2 + 1)$ we obtain a substantial power saving. In fact, even the rather generic (though highly non-trivial) Jacquet-Shalika bounds $N_v(\sigma, \mathcal{F}_I(q)) = 0$ for $\sigma > 1/2$ \cite{JS} show for $n \geq 3$ that cuspidal representations are always fairly far away from the trivial representation. Unfortunately it is not clear how to combine the Luo-Rudnick-Sarnak approach with the  present techniques. 
 
 The theorem remains true for $n=2$ (by a slightly different proof), but it is known in this case (see \cite{Iw} for $v=\infty$, and the proof for finite $v$ is similar) and recovers Selberg's $3/16$ bound for exceptional eigenvalues. For $n=3$ and $v=\infty$ this is \cite[Theorem 4]{BBM}. As mentioned above, for larger $n$ Theorem \ref{thm1} is completely new. As we shall outline below it appears to be the limit of what is available by any trace formula approach, even in the case  $n=2$ nothing better is known. 
 
  The proof is based on a careful analysis of the arithmetic side of the Kuznetsov formula with a test function on the spectral side that blows up on exceptional Langlands parameters at $v$ (and therefore increases the complexity on the arithmetic side). We denote by $\lambda_{\pi}(m)$ the $m$-th Hecke eigenvalue of $\pi \in \mathcal{F}_I(q)$.  
 
 \begin{theorem}\label{thm2} Keep the assumptions and notation of Theorem \ref{thm1}. Let $m\in\Bbb{N}$ be coprime to $q$ and $Z \geq 1$. Then
 \begin{equation*}%\label{2a}
 \sum_{\pi \in \mathcal{F}_I(q)} |\lambda_{\pi}(m)|^2 Z^{2\sigma_{\pi}(\infty)}  \ll_{I, n, \varepsilon} q^{n-1+\varepsilon}
 \end{equation*}
% and
%\begin{equation}\label{2b}
%\sum_{\pi \in \mathcal{F}_I(q)} Z^{2\sigma_{\pi}(\infty)} \ll_{I, n, \varepsilon} q^{n-1+\varepsilon}
%\end{equation}
 uniformly in $mZ \ll q^2$ for a sufficiently small implied constant (depending on $I$ and $n$). 
 \end{theorem}
 
We shall see in Lemma \ref{Lenstra} below that $|\lambda_{\pi}(p^{\nu})|^2$ is often as big as $p^{2\nu \sigma_{\pi}(p)}$ for a prime $p$ and $\nu \in \Bbb{N}$, so that the ``test function'' $|\lambda_{\pi}(m)|^2 Z^{2\sigma_{\pi}(\infty)} $ treats finite  places  and the infinite place essentially on equal footing. %It gives finer information on the spectral side, but  blows up the length of the arithmetic side of the Kuznetsov formula, cf.\ \eqref{heur}. 

Let us roughly sketch how one may hope to arrive at Theorem \ref{thm2}. %\eqref{2a} (the bound \eqref{2b} is  a non-archimedean analogue). 
 Since the Laplacian eigenvalue is fixed, the Whittaker transforms in the Kuznetsov formula play no major role, and the battle is decided on the level of Kloosterman sums. Very roughly, the Kuznetsov formula takes the shape
 \begin{equation}\label{heur}
 \begin{split}
& \frac{1}{|  \mathcal{F}_I(q)| } \sum_{\pi \in \mathcal{F}_I(q)} |\lambda_{\pi}(m)|^2Z^{2\sigma_{\pi}(\infty)}  \,\,\, ``\approx"  \,\,\, 1 + \sum_{\text{id} \not= w \in W} \sum_{\substack{q \mid c_1, \ldots, q \mid c_{n-1}\\ c_1, \ldots, c_{n-1} \ll mZ}} \frac{S_{q, w}(M, M, c)}{c_1 \cdots c_{n-1}}
 \end{split}
 \end{equation}
 where $W$ is the Weyl group of permutation matrices, $M = (m, 1, \ldots, 1) \in \Bbb{Z}^{n-1}$  and $S_{q, w}(M, M, c)$ is a certain generalized Kloosterman sum, defined in \eqref{klo} below, associated with the Weyl element $w$ and moduli $c = (c_1, \ldots, c_{n-1})$. If $mZ \ll q$, then the off-diagonal term vanishes completely and we are done. We will use this observation in Theorem \ref{thm4} below. This range of $mZ$  recovers the ``convexity bound'' with the value \eqref{c}.  For larger values of $mZ$ and  stronger density results we must deal with the Kloosterman sums appearing in the off-diagonal term and  improve on the trivial bound $|S_{q, w}(\ast, \ast, c)| \ll (c_1 \cdots c_{n-1})^{1+\varepsilon}$, see \eqref{triv} below. % It is the possibility of using non-trivial bounds for Kloosterman sums (if they are available) that opens the door to ``subconvex'' density results. 
  To obtain such  bounds for general groups is a famous open problem. In an ideal world we would have complete Weil-type square root cancellation $|S_{q, w}(\ast, \ast, c)| \ll (c_1 \cdots c_{n-1})^{1/2+\varepsilon}$ (at least under certain coprimality assumptions) which allows us to take $mZ$ as large as $q^2$. This square-root cancellation implies the statement of Theorem \ref{thm2} and   a density hypothesis halfway between the trivial representation and the Ramanujan conjecture. (Additional square root cancellation in the $c_1, \ldots, c_{n-1}$ sum would give the full Ramanujan conjecture, but this is of course not in the cards).
 
  Interestingly, this heuristic sketch turns out to be quite far from the truth.  Square root cancellation for the size of Kloosterman sums may  fail badly, so that we have to arrive at Theorem \ref{thm2} by a rather different analysis. The key lemma is the following, which seems to be the first explicit analysis of general ${\rm GL}(n)$ Kloosterman sums beyond hyper-Kloosterman sums \cite{Fr} associated to the Weyl element $w = \left(\begin{smallmatrix} &  1\\  I_{n-1} &  \end{smallmatrix}\right)$, where $I_{n}$ denotes the $n$-by-$n$ identity matrix.  
 
 \begin{theorem}\label{thm3} Let $q$ be a prime and let $M, N \in \Bbb{Z}^{n-1}$ with entries  coprime to $q$ (in particular non-zero). Let $n \geq 3$ and let $w \in W$. Then $S_{q, w}(M, N, (q, \ldots, q)) = 0$ unless   
 \begin{equation}\label{wast}w = w_{\ast} :=  \left(\begin{smallmatrix} & & 1\\ & I_{n-2} &\\ 1 &&  \end{smallmatrix}\right) 
 \end{equation} in which case 
$S_{q, w}(M, N, (q, \ldots, q)) = q^{n-2}.$
 \end{theorem}
 
Note that this is in sharp contrast to the case $n=2$, where a Kloosterman sum to prime modulus $q$ has no closed evaluation. The key point is that by multiplicativity the Kloosterman sums in \eqref{heur} contain $S_{q, w}(\ast, \ast, (q, \ldots, q))$ as a large chunk. The critical case is the term corresponding to $w = w_{\ast}$ where the Kloosterman sum is much bigger than the product of the square root of the moduli  (if $n > 3$). Luckily in this case the remaining
 piece with moduli $(c_1/q, \ldots, c_{n-1}/q)$ comes with additional savings since the Weyl element $w_{\ast}$ imposes certain relations among the $c_j$.   That the critical Weyl element is not the long Weyl element, but  rather the  permutation $1 \leftrightarrow n$ that is relatively ``close'' to the identity, may also be quite surprising in this context. This analysis is sensitive to $q$ being prime. While the general technique can be applied in a rather broad context (see e.g. \cite{AB, Ma}), the estimation in Theorem \ref{thm3} is particularly designed to the particular setup of $\Gamma_0(q)$, $q$ prime. 
 
Theorem \ref{thm2} and variations of it have other applications of which we mention here only one, namely a large sieve inequality. 
\begin{theorem}\label{thm4} Let $q$ be prime and  $(\alpha(m))$  any sequence of complex numbers. Then
$$\sum_{\pi \in \mathcal{F}_I(q)} \Big| \sum_{\substack{m \leq x\\ (m, q) = 1}} \alpha(m) \lambda_{\pi}(m) \Big|^2 \ll_{I, n, \varepsilon} q^{n-1+\varepsilon} \sum_{\substack{m \leq x\\ (m , q) = 1}} |\alpha(m)|^2$$
uniformly in $x \ll q$ for a sufficiently small implied constant  (in terms of $I$ and $n$). 
\end{theorem}

This result holds (with literally the same proof) for all $q \in \Bbb{N}$. For comparison, Venkatesh \cite[Theorem 1]{Ve} obtained this with $x \leq q^{1/(2n-2)}$. A simple corollary is the following best-possible bound for a second moment of $L$-functions on the critical line:
\begin{cor}\label{cor5} For $q$ prime and $t \in \Bbb{R}$ we have
$$\sum_{\pi \in \mathcal{F}_I(q)} |L(1/2 + it, \pi)|^2 \ll_{I, t, n, \varepsilon} q^{n-1 + \varepsilon}.$$
\end{cor}

The author would like to thank Farrell Brumley for encouragement and numerous discussions on the subject. 
 
\section{Basic notation} Let $U \subseteq {\rm GL}_n$ be the subgroup of unipotent upper triangular matrices. The Haar measure on $U(\Bbb{R})$ is given by    $\dd x = \prod_{1 \leq i < j \leq n} \dd x_{ij}$. As before let $W$ be the Weyl group of permutation matrices; we identify a permutation matrix $w = (w_{ij}) \in W$ with the  permutation $ i \mapsto j$ for $w_{ij} = 1$. 
For $w \in W$ we define $$U_w = w^{-1} U^{\top} w \cap U$$ 
and we continue to write    $\dd x$  for the induced measure on the subgroup $U_{w}(\Bbb{R})$. (Here $^{\top}$ denotes the transpose, so that $U^{\top}$ is the set of unipotent lower triangular matrices.) 

As 
\begin{equation}\label{ww}
w^{-1} (x_{ij}) w = (x_{w^{-1}(i), w^{-1}(j)}),
\end{equation}
 the group $U_w$ has entries at $(i, j)$ with $i< j$ exactly when $w^{-1}(i)  > w^{-1}(j)$ (since $U^{\top}$ consists of lower triangular matrices). Equivalently, $U_w$ has entries at 
\begin{equation}\label{entries}
   (w(i), w(j)) \text{ for } 1 \leq j < i \leq n  \text{ whenever } w(i) < w(j).
   \end{equation} Let $V \subseteq {\rm GL}_n$ be the group of   diagonal matrices with entries $\pm 1$.  

For  $N   \in \Bbb{Z}^{n-1}$ we define a character $\theta_N : U(\Bbb{R})/U(\Bbb{Z}) \rightarrow S^1$ by 
\begin{equation}\label{character}
   \theta_N(x) = e(N_{n-1} x_{12} + \ldots + N_{1} x_{n-1, n})% = \theta_{(1, \ldots, 1)}(N x N^{-1})
.
\end{equation}
 For $v \in V$ we write $\theta_N^v(x) = \theta_N(v^{-1} x v)$ (note that $v^{-1} U v = U$).  If $N = (1, \ldots, 1)$, we drop it from the notation of the character. 

%For $q \in \Bbb{N}$ we write $\Gamma_0(q) \subseteq {\rm SL}_n(\Bbb{Z})$ for subgroup with entries $x_{n, 1}, \ldots, x_{n, n-1}$ divisible by $q$. 

Let $T \subseteq {\rm GL}_n$ be the diagonal torus. We embed   $y = (y_1, \ldots, y_{n-1}) \in \Bbb{G}_m^{n-1}$ into $T$ as 
\begin{equation}\label{iota}
\iota(y) = \text{diag}(y_{n-1} \cdots  y_1, \ldots, y_2y_1, y_1, 1).  %= (y_1 \cdots y_{n-j})_{1 \leq j \leq n}
\end{equation}
We multiply two elements in $y, y' \in \Bbb{G}_m^{n-1}$ componentwise, written $y \cdot y'$, so that $\iota$ is a homomorphism. 
  We denote the image of $\Bbb{R}_{>0}^{n-1}$ in $T$ by $\tilde{T}(\Bbb{R})$. Then $\mathcal{H} = U(\Bbb{R}) \tilde{T}(\Bbb{R})$ is the generalized upper half plane in the sense of \cite[Chapter 1]{Go}. We identify $\mathcal{H}$ with ${\rm GL}_n(\Bbb{R})/{\rm O}_n(\Bbb{R}) {\rm Z}^{+}$ where ${\rm Z}^{+} \cong  \Bbb{R}_{>0}$ is the subgroup of diagonal scalar matrices  with positive entries. For $g = x y k \alpha \in {\rm GL}_n(\Bbb{R})$ with $x \in U(\Bbb{R})$, $y \in \tilde{T}(\Bbb{R})$, $k \in {\rm O}_n(\Bbb{R})$, $\alpha \in {\rm Z}^+$, we write ${\rm y}(g) = \iota^{-1}y \in \Bbb{R}_{>0}^{n-1}$ for $(n-1)$-tuple of  Iwasawa $y$-coordinates. In particular, for $g = \text{diag}(y_1, \ldots, y_n)$ with positive $y_j$ we have 
\begin{equation}\label{y}
  {\rm y}(g) = \Big(\frac{y_{n-1}}{y_n}, \ldots, \frac{y_1}{y_2}\Big) \in \Bbb{R}_{>0}^{n-1}.
  \end{equation}
 For $w \in W$, $y \in \Bbb{R}_{>0}^{n-1}$ we write $${\rm y}(w\iota(y)^{-1} w^{-1}) = {}^wy  = ({}^wy_1, \ldots, {}^wy_{n-1})$$ for the Iwasawa $y$-coordinates of $w\iota(y)^{-1} w^{-1}$. Explicitly, combining \eqref{ww} with $w^{-1}$ in place of $w$, \eqref{iota} and \eqref{y}, we have
 \begin{equation}\label{yNw}
{}^wy = \Big( \frac{y_1 \cdots y_{n - w(n-j+1)}}{y_1 \cdots y_{n- w(n-j)}}\Big)_{1 \leq j \leq n-1}.
 \end{equation}
  
 %In particular, every $g \in {\rm GL}_n(\Bbb{R})$ has a unique representative modulo ${\rm O}_n(\Bbb{R}) \Bbb{R}_{>0}$ in $\mathcal{H}$.

For $\alpha\in \Bbb{C}^{n-1}$, $y \in \Bbb{R}_{>0}^{n-1}$ we write $y^{\alpha} = y_1^{\alpha_1} \cdots y_{n-1}^{\alpha_{n-1}}\in \Bbb{C}$. 
Let 
\begin{equation}\label{eta}
\eta = (\eta_1, \ldots, \eta_{n-1}) =   \Big(\frac{1}{2}j(n-j)\Big)_{1 \leq j \leq n-1}. 
\end{equation}
We define a measure on $\Bbb{R}_{>0}^{n-1}$   by $\dd^{\ast}y = y^{-2\eta} \frac{dy_1}{y_1} \cdots \frac{dy_{n-1}}{y_{n-1}}$ and correspondingly an inner product   by $$\langle f, g \rangle = \int_{\Bbb{R}_{>0}^{n-1}} f(y) \bar{g}(y) \dd^{\ast}y.$$  We denote the push forward of  $\dd^{\ast}y$ to $\tilde{T}(\Bbb{R})$ by $\iota$ also by $\dd^{\ast}y$. 
Then $\dd x \, \dd^{\ast}y$ is a left ${\rm GL}_n(\Bbb{R})$ invariant measure on $\mathcal{H}$. 

 We define a different embedding of $\Bbb{R}_{>0}^{n-1}$ into $T(\Bbb{R})$  by $$c = (c_1, \ldots, c_{n-1}) \mapsto c^{\ast} = \text{diag}(1/c_{n-1}, c_{n-1}/c_{n-2}, \ldots, c_2/c_1, c_1).$$
From \eqref{y},  it is useful to observe that   
 \begin{equation}\label{c-y}
   {\rm y}(c^{\ast})  = \Big(\frac{c_{j-1}c_{j+1}}{c_j^2}\Big)_{1 \leq j \leq n-1}
   \end{equation}
 where $c_0 = c_n = 1$, and a simple computation shows
 \begin{equation}\label{ceta}
   {\rm y}(c^{\ast})^{\eta} = (c_1 \cdots c_{n-1})^{-1}.
  \end{equation}

\section{Auxiliary results}

As the Iwasawa decomposition (with the compact group on the right) is the Gram-Schmidt orthogonalization of rows starting with the last row, we can compute ${\rm y}(g)$  explicitly.  For $1 \leq j \leq n$ let $\Delta_j = \Delta_j(g)$ be the volume of the parallelepiped spanned by last $j$ rows  of $g$.  Then 
  $$g \equiv \left(\begin{smallmatrix} \Delta_n/\Delta_{n-1} & \ast &\cdots  & \ast \\ & \Delta_{n-1}/\Delta_{n-2} & \cdots  & \ast \\ & & \ddots & \vdots \\ &&&\Delta_1\end{smallmatrix}\right) \quad (\text{mod } O_n(\Bbb{R})),$$
so that by \eqref{y} we have 
\begin{equation}\label{gram}
  {\rm y} (g) = \Bigl(\frac{\Delta_{j+1}(g) \Delta_{j-1}(g)}{\Delta_{j}(g)^2}\Big)_{1 \leq j \leq n-1}
\end{equation}  
with  the convention $\Delta_{0}(g) = 1$. By \cite[Corollary 4.2 and p.\ 11]{FP} (or by hand) we confirm the inversion formula for the following $(n-1)$-by-$(n-1)$ tridiagonal Toeplitz matrix
\begin{equation}\label{inverse}
   \left( \begin{smallmatrix} -2 & 1 &  &\\ 1 & -2 & 1&  \\   & \ddots & \ddots& \ddots\\ & & 1 &-2 &1\\ & & & 1 & -2 \end{smallmatrix}\right)^{-1} = \Big(-s(i, j)\Big)_{ij}, \quad s(i, j) = \frac{1}{n} \begin{cases} i(n-j), & i \leq j, \\ j(n-i), & i > j.  \end{cases} 
  \end{equation}%\min\big(i(n-j), (n-i)j\big)$$
%(see \cite[Corollary 4.2 and p.\ 11]{FP} or check quickly by hand). 
%$$s(i, j) = \frac{1}{n} \min\big(ij, (n-i)(n-j)\big).$$
Therefore, given ${\rm y}(g) = (Y_1, \ldots, Y_{n-1})$ and $\Delta_n(g) = |\det(g)|$ we can solve \eqref{gram} explicitly for $\Delta_1, \ldots, \Delta_{n-1} > 0$   getting
\begin{equation}\label{delta}
\Delta_j(g) = |\det(g)|^{j/n} \prod_{i=1}^{n-1} Y_i^{-s(i, j)}. 
\end{equation}
%which follows from the readily checked identity
%$$s(i, n-j+1) + s(i, n-j-1) - 2 s(i, n-j) + \delta_{i=j} = 0.$$ 

Our first lemma will be used to bound the moduli $c$ on the arithmetic side of the Kuznetsov formula. 

\begin{lemma}\label{volume} Let   $w\in W$, $x \in U_w(\Bbb{R})$, $y, c, B \in \Bbb{R}_{>0}^{n-1}$. Write ${\rm y}\big(\iota(B)c^{\ast}w x\iota(y) \big) = Y \in \Bbb{R}_{>0}^{n-1}$ and $A = \iota(B) c^{\ast}$.  Then
\begin{equation*}%\label{first}
\begin{split}
&c_j  \ll_{y, Y} %c_j \Delta_j(wx) = 
\prod_{i=1}^{n-1}  B_i^{s(i,j)}\quad \text{and} \quad 
%\end{equation*}
%$$c_j \leq c_j \Delta_j(wx) \ll X^{ s(n-1, n-j) \frac{n(n+1)}{2}-\sum_{i=1}^{n-1} s(i, n-j) }  \prod_{i=1}^{n-1} B_i^{s(i, n-j)}$$
%for $1 \leq j \leq n-1$ and  $\Delta_n(wx) = 1$. 
%and
%\begin{equation*}%\label{second}
 1 \leq  \Delta_j(wx)  \ll_{y, Y} \prod_{i=1}^{n-1} {\rm y}(A)_i^{s(i,j)} 
\end{split}
\end{equation*}
for $1 \leq j \leq n-1$. 
\end{lemma}

 \textbf{Proof.}  We have %$\Delta_j(w x) \geq 1$ since one of the minors is always 1, so 
 \begin{equation}\label{1}
 \Delta_j\big(\iota(B)c^{\ast}w x\iota(y) \big)  = \Delta_j(wx \iota(y) ) c_j   \prod_{i=1}^j  B_1 \cdots B_{i-1}
 \end{equation} %\geq c_j \prod_{i=1}^j y_1 \cdots y_{n- w(n+1-i)} B_1 \cdots B_{i-1} $$
since the diagonal matrix $\iota(B)c^{\ast}$ multiplies the rows of $w x\iota(y)$ by the corresponding diagonal entries. 
Clearly $\Delta_j(wx) \geq 1$  since one of the minors is always 1, and clearly  
 \begin{equation}\label{35}
 |\det(\iota(B)c^{\ast}w x\iota(y) )| = |\det(\iota(B) \iota(y))| \asymp_y B_1^{n-1} \cdots B_{n-2}^2 B_{n-1}.
 \end{equation}
  From \eqref{1}, and from \eqref{delta} in combination with \eqref{35}, we therefore obtain
 \begin{equation}\label{long}
 \begin{split}
 c_j  \leq c_j \Delta_j(wx)& \asymp_y  c_j \Delta_j(wx\iota(y)) =  \Delta_j\big(\iota(B)c^{\ast}w x\iota(y) \big)    \prod_{i=1}^{j-1} B_i^{-(j-i)}\\
 &  \asymp_{Y, y}   \prod_{i=1}^{n-1}  B_i ^{(n-i)j/n}  \prod_{i=1}^{j-1} B_i^{-(j-i)}  = \prod_{i=1}^{n-1} B_i^{s(i, j)}. % = \prod_{i=1}^{n-1} (B_i/Y_i)^{s(i, j)}
 \end{split}
 \end{equation}
 This shows the first statement of the lemma, and the proof of the second is completed by observing that  \eqref{c-y} and \eqref{inverse} imply $$\frac{1}{c_j} = \prod_{i=1}^{n-1} {\rm y}(c^{\ast})_i^{s(i, j)}$$
for $1 \leq j \leq n-1$, so that \eqref{long} implies $ \Delta_j(wx)  \asymp_{y, Y}  \prod_{i=1}^{n-1} ({\rm y}(c^{\ast})_iB_i)^{s(i, j)}$ as desired.  \hfill $\square$\\
% with the Kronecker symbol, and by \eqref{wNy} also
% $$\prod_{i=1}^{n-1} {}^wy_i^{-s(i, j)} = \prod_{i = 1}^n (y_1 \cdots y_{n-w(i)})^{s

We shall see in a moment that the only Weyl elements contributing to the Kuznetsov formula are of the form
\begin{equation}\label{w}
w = \left( \begin{matrix} &  & & I_{d_1}\\  & & I_{d_2} &\\ & \Ddots & &\\ I_{d_r} & && \end{matrix}\right)
\end{equation}
with identity matrices $I_{d_j}$ of dimension $d_j$ (i.e.\ $d_1 + \ldots  +d_r = n$), 
so without loss of generality we restrict our attention to such matrices. The following technical result computes the Jacobi determinant for a certain change of variables. 

\begin{lemma}\label{lemma1} Let $N \in \Bbb{N}^{n-1}$, $w\in W$ of the form \eqref{w}. For $x \in U_w(\Bbb{R})$ define $x' = \iota(N) x \iota(N)^{-1} \in U_w(\Bbb{R})$. Then
$$\frac{\dd x'}{\dd x} = ({}^wN)^{\eta} N^{\eta}$$
where the left hand side denotes the Jacobi determinant $\det Dx'(x)$. 
\end{lemma}

%\textbf{Remark.} The result holds for arbitrary $w$, but we do 
 
\textbf{Proof.} Since $\iota(N)(x_{ij}) \iota(N)^{-1} = (x_{ij} N_1 \cdots N_{n-i}(N_1 \cdots N_{n-j})^{-1})_{ij}$ and recalling \eqref{entries} and \eqref{yNw}, we have to show 
\begin{equation}\label{toshow}
 \prod_{\substack{1 \leq j < i \leq n\\ w(i) < w(j)}} N_{n- w(j) + 1} \cdots N_{n - w(i)} = \prod_{j=1}^{n-1}  \Big( N_j\frac{N_1 \cdots N_{n - w(n-j+1)}}{N_1 \cdots N_{n- w(n-j)}}\Big)^{\eta_j}
 \end{equation}
 for an arbitrary $w$ as in \eqref{w}. We use induction on $r$ and write $w' = \left(\begin{matrix} & w\\ I_d & \end{matrix}\right)$, so that $n+d - w'(j) = n - w(j)$ for all $1 \leq j \leq n$.  We call $L(w)$ the left hand side of \eqref{toshow} and $R(w)$ the right hand side. We consider first the quotient $L(w')/L(w)$. The pairs $1 \leq j < i \leq n$ cancel, and for $i > n$ only $j \leq n$ satisfy the summation condition $w'(i) < w'(j)$. We conclude
%Since $w$ is of the form \eqref{w}, we may call the the left hand side $L(I_{d_1}, \ldots, I_{d_r})$ and the right hand side $R(I_{d_1}, \ldots, I_{d_r})$. We use induction by $r$. %We have trivially $L(I_{d_1}) = R(I_{d_1}) = 1$. Suppose that \eqref{toshow} holds for some $r$ and write $d_1 + \ldots + d_r = n$ and $d = d_{r+1}$. Then
%If we write $d_1 + \ldots + d_r = n$ and $d = d_{r+1}$, then
\begin{equation}\label{LHS}
\frac{L(w')}{L(w)}  =  \prod_{j=1}^n \prod_{i=n+1}^{n+d} N_{n+d- w'(j) + 1} \cdots N_{n +d- w'(i)} = \prod_{j=1}^{n-1} N_j^{dj} \prod_{j=n}^{n+d-1} N_j^{n(n+d-j)}. 
\end{equation}
On the other hand,  
$$R(w) = \prod_{j=1}^{n-1}N_j^{\eta_j} \prod_{i = 1}^n (N_1 \cdots N_{n- w(i)})^{\eta_{n-i+1} - \eta_{n-i}} = \prod_{j=1}^{n-1}N_j^{\frac{j(n-j)}{2}} \times \prod_{i = 1}^n (N_1 \cdots N_{n- w(i)})^{\frac{2i-n-1}{2}},$$
so %again writing $d_1 + \ldots + d_r = n$ and $d = d_{r+1}$ we see that 
$R(w')/
R(w)$ equals
\begin{displaymath}
\begin{split}
%&\frac{R(I_{d_1}, \ldots, I_{d_{r+1}})}{R(I_{d_1}, \ldots, I_{d_r})} \\
&  \prod_{j=1}^{n-1} N_j^{\frac{j(n+d - j)}{2} - \frac{j(n-j)}{2}} \prod_{j=n}^{n+d-1} N_j^{\frac{j(n+d - j)}{2} }\times \prod_{i = 1}^n (N_1 \cdots N_{n- w(i)})^{-\frac{d}{2}}\prod_{i = n+1}^{n+d} (N_1 \cdots N_{n+d- w'(i)})^{\frac{2i-n-d-1}{2}}\\
&= \prod_{j=1}^{n-1} N_j^{\frac{dj}{2}}  \prod_{j=n}^{n+d-1} N_j^{\frac{j(n+d - j)}{2} }\times  \prod_{j=1}^{n-1} N_j^{-\frac{(n-j)d}{2}}  \prod_{j=1}^{n-1} N_j^{\sum_{i=n+1}^{n+d} \frac{2i-n-d-1}{2}}  \prod_{j=n}^{n+d-1} N_j^ {\sum_{i=n+1}^{2n+d-j} \frac{2i-n-d-1}{2}}
\end{split}
\end{displaymath}
which is easily seen to equal the right hand side of \eqref{LHS}. Since trivially $L(I_{d}) = R(I_{d}) = 1$, the induction is complete.  \hfill $\square$

%$$\frac{\dd x'}{\dd x}  = \prod_{\substack{1 \leq j < i \leq n\\ w(i) < w(j)}} N_{n- w(j) + 1} \cdots N_{n - w(i)}.$$ \\ 
%We \eqref{yNw} we have to show that this equals
%$$  \prod_{j=1}^{n-1}  \Big( N_j\frac{N_1 \cdots N_{n - w(n-j+1)}}{N_1 \cdots N_{n- w(n-j)}}\Big)^{\eta_j}$$

\begin{lemma}\label{int-lemma} Let $B \in \Bbb{R}_{>0}^{n-1}$, $w = w_{\ast}\in W$ as in \eqref{wast}. Then  
$$\text{{\rm vol}}\{x \in U_w(\Bbb{R}) \mid \Delta_j(wx) \leq B_j, 1 \leq j \leq n-1\} \ll_{\varepsilon} (B_1 \cdots B_{n-1})^{1+\varepsilon}$$
for any $\varepsilon > 0$. 
\end{lemma}

\textbf{Proof.} We can assume without loss of generality that $B_j \geq 1$, otherwise the volume is 0 as seen in the proof of Lemma \ref{volume}. For $x \in U_w(\Bbb{R})$ we have 
$$wx = \left(\begin{matrix} &&& 1\\ &&& x_{2, n} \\&& I_{n-2}& \vdots\\ &   & & x_{n-1, n}
\\1 & x_{12} & \ldots & x_{1, n}\end{matrix}\right),$$
so that by considering the lower right minors we obtain in particular the inequalities 
$$|x_{1 n}| \leq B_1, \quad  \Big|-x_{1, n} +  \sum_{i = n+1-j}^{n-1} x_{1i} x_{i,n} \Big| \leq B_{j}, \quad j = 2, \ldots, n-1,$$
and we also have $|x_{ij}| \leq b := 1+ \max(B_1, \ldots, B_{n-1})$. If $I \subseteq \Bbb{R}$ is any interval of length $|I| \geq 1$, then  
\begin{displaymath}
\begin{split} \text{vol}\{(x, y) \in [-b, b]^2 : xy \in I\}& \leq   \int_{|y| \leq b} \min\Big(\frac{|I|}{|y|} , 2b\Big)\dd y \leq 4|I| + \int_{|I|/b \leq |y| \leq b} \frac{|I|}{|y|} \dd y\\
& \leq  4 |I| (1 + \log b ). 
\end{split}
\end{displaymath}
Thus if $|x_{1n}| \leq B_1$ is fixed, the volume of $(x_{1, n-1}, x_{n-1, n})$ is $O(B_2 \log b)$, and if these are fixed,   the volume of $(x_{1, n-2}, x_{n-2, n})$ is $O(B_3 \log b)$, etc. Inductively we obtain the desired bound.  \hfill $\square$\\

Most likely the statement holds for all $w$, but the proof is particularly simple for $w_{\ast}$ which is all we need.

\section{Kloosterman  sums} 
Properties of Kloosterman sums for ${\rm SL}_n(\Bbb{Z})$ have been obtained and summarized in \cite{Fr}. They generalize in an obvious way to the congruence subgroup $\Gamma_0(q)$. 
 The Bruhat decomposition gives ${\rm GL}_n( \Bbb{Q}) = \bigcup_{w\in W} G_w(\Bbb{Q})$ with $G_w := U  T  w U_w $ as a disjoint union. 
Let $N, M, c \in \Bbb{Z}^{n-1}$, $w \in W$, $v \in V$. Then  provided that
\begin{equation}\label{compat}
  \theta_M(c^{\ast} w x w^{-1} (c^{\ast})^{-1}) = \theta_N^v(x)
\end{equation}
for all $x \in w^{-1}U(\Bbb{Q}) w \cap U(\Bbb{Q})$ [this set is a ``complement'' of $U_w$ in $U$], the   Kloosterman sum
\begin{equation}\label{klo}
S^v_{q, w}(M, N, c) = \sum_{x  c^{\ast} w y \in U(\Bbb{Z})\backslash G_w(\Bbb{Q}) \cap \Gamma_0(q)/U_w(\Bbb{Z}) } \theta_M(x )\theta_N^v(y)
\end{equation}
is well-defined, see \cite[Proposition 1.3]{Fr}.  If \eqref{compat} is not met, we define $S^v_{q, w}(M, N, c) = 0$. If $v = \text{id}$, we drop it from the notation. By \cite[p.\ 175]{Fr}, the Kloosterman sum is non-zero only if $w$ is of the form \eqref{w}. 
%\begin{equation}\label{w}
%w = \left( \begin{matrix} &  & & I_{d_1}\\  & & I_{d_2} &\\ & \Ddots & &\\ I_{d_r} & && \end{matrix}\right)
%\end{equation}
 If $\gamma = x_1 c^{\ast} w x_2\in \Gamma_0(q)$ is a matrix occurring in the sum on the right hand side of \eqref{klo}, then any minor of $\gamma$, and hence of $c^{\ast} w$, obtained by deleting at least the first row and the last column is divisible by $q$. Hence if $w$ is of the form \eqref{w}, then the summation condition in \eqref{klo} can only be met if
\begin{equation}\label{divi}
 q \mid c_1, \quad q\mid c_2, \quad \ldots, \quad q \mid c_{n -d_1}.
\end{equation}
%We write $\bar{N} = (N_{n-1}, \ldots, N_1)$ and 
Observing that  
\begin{equation}\label{char}
  \theta_{M}(x) = \theta(\iota(M) x \iota(M)^{-1})
  \end{equation}
and recalling \eqref{c-y}, we see that   \eqref{compat} is equivalent to
\begin{equation}\label{comp-equiv}
 M_{n-i} \frac{c_{n-i+1} c_{n-i-1}}{c_{n-i}^2} =  \frac{v_{w(i) + 1}}{v_{w(i)}} N_{n-w(i)}
 \end{equation}
for all   $1 \leq i \leq n-1$ satisfying   $w(i) + 1 = w(i+1)$ with the above convention $c_0 = c_n = 1$ and $v = \text{diag}(v_1, \ldots, v_n)$.  If $w$ is of the form \eqref{w}, these are precisely the $i \not\in \{d_1, d_1 + d_2, \ldots, d_1 + d_2 + \ldots + d_{r-1}\}$. If $w = \text{id}$, then $xc^{\ast} w y = xc^{\ast} y$ can only be in $\Gamma_0(q)$ if $c_1 = \ldots = c_{n-1} = 1$, in which case we conclude from \eqref{comp-equiv} that $M_j = \pm N_j$. 

Kloosterman sums for ${\rm SL}_n(\Bbb{Z})$ enjoy certain  multiplicativity properties in the moduli, cf.\  \cite[Proposition 2.4]{Fr}. We state only one particular case. Let $q$ be  prime, suppose that $(c_1\cdots c_{n-1}, q) = 1$ and write $qc = (qc_1, \ldots, qc_{n-1})$. Suppose that $w(1) = n$ and $w(n) = 1$. Then 
%in the following sense: if $(c_1'\cdots c_{n-1}', c_1''\cdots c_{n-1}'') = 1$ and $c = (c_1'c_1'', \ldots, c_{n-1}'c_{n-1}'')$, then by \cite[Proposition 2.4 (and Proposition 2.3)]{Fr} we have
\begin{equation}\label{mult}
  S^v_{q, w}(M, N, qc) = S^v_{q, w}(M, N' , (q, \ldots, q)) S^v_{1, w}(M, (\bar{q}N_1, N_2, \ldots, N_{n-2}, \bar{q}N_{n-1}) , c)
\end{equation}
with
\begin{displaymath}
\begin{split}
& N_{n-i}' \equiv N_{n-i} c_{n - w(i)} c_{n - w(i+1) + 1} \overline{c_{n - w(i)+1} c_{n - w(i+1) } } \, (\text{mod } q).\\
%& N_{n-i}'' \equiv N_{n-i} c'_{n - w(i)} c'_{n - w(i+1) + 1} \overline{c'_{n - w(i)+1} c'_{n - w(i+1) } } \, (\text{mod } q).
\end{split}
\end{displaymath}
 
%We spell out two instances of \eqref{compat} explicitly. Suppose that $w$ is of the form \eqref{w} with $\dim I_1 = \ell > 1$. Then the  entry $x_{n+1 - \ell, n}$ of $x \in w^{-1} U w \cap U$ is not restricted, and this becomes the entry with indices $(1, 2)$ of $wxw^{-1}$. Comparing coefficients, we conclude
%\begin{equation}\label{compat1}
%\frac{N_{n-i} c_{n-2}}{c_{n-1}^2} = M_{n+1-\ell}.
%\end{equation}
%Similarly, $w$ is of the form \eqref{w} with $\dim I_r = \ell > 1$ we obtain
%\begin{equation}\label{compat1}
%\frac{N_{n-1} c_{ 2}}{c_{ 1}^2} = M_{\ell-1}.
%\end{equation}
By \cite[Theorem 0.3(i)]{DR} we have the trivial bound
\begin{equation}\label{triv}
|S^v_{q, w}(M, N, c) | \leq    | U(\Bbb{Z})\backslash G_w(\Bbb{Q}) \cap {\rm SL}_n( \Bbb{Z})/U_w(\Bbb{Z})|  \ll (c_1 \cdot \ldots \cdot c_{n-1})^{1+\varepsilon}. 
\end{equation}\\

We now give the \textbf{proof of  Theorem \ref{thm3}} from the introduction, which  is the first non-trivial bound for a ${\rm GL}_n$ Kloosterman sum other than a hyper-Kloosterman sum.  For $n=3$, the  statement is essentially contained in \cite[Lemma 6(c)]{BBM}. We wish to compute $S_{q, w}(M, N, (q, \ldots, q)) = 0$ where $M, N \in \Bbb{Z}^{n-1}$ have entries coprime to $q$.

As mentioned before, we can assume that $w$ is of the form \eqref{w}, otherwise the Klooster\-man sum vanishes by definition. Next assume that $d_1 > 1$ in \eqref{w}. Applying \eqref{comp-equiv} with $i=1$ we obtain $M_{n-1} = \pm N_{d_1 - 1} q$, a contradiction. In the same  way we exclude the case $d_r > 1$. For $w$ of the form \eqref{w} with $d_1 = d_r = 1$ and $c^{\ast} = \text{diag}(1/q, 1, \ldots, 1, q)$, we recall the definition \eqref{klo} and consider $\gamma = x  c^{\ast} w y \in G_w(\Bbb{Q}) \cap \Gamma_0(q)$ with uniquely determined $x  \in U(\Bbb{Z})\backslash U(\Bbb{Q})$, $y \in U_w(\Bbb{Q})/U_w(\Bbb{Z})$. A system of representatives for $U(\Bbb{Z})\backslash U(\Bbb{Q})$ consists of matrices with rational entries in $[0, 1)$ above the diagonal, and similarly we choose a system of representatives 
of $U_w(\Bbb{Q})/U_w(\Bbb{Z})$ where all relevant entries are restricted to $[0, 1)$. We now determine those representatives $x, y$  that satisfy $xc^{\ast} w y \in \Gamma_0(q)$. We have 
$$x = \left(\begin{matrix}1 & \ast & \ast &\cdots& \ast & x_{1, n}\\ & 1 & \ast&\cdots &\ast &x_{2,n}\\ & & \ddots& && \vdots\\ &&&1  &\ast& x_{n-2, n}\\ & & && 1& x_{n-1, n} \\ && & & &1  \end{matrix}\right) , \quad\quad c^{\ast} w y =  \left( \begin{matrix} &&  && & 1/q\\ & && & I_{d_2} & y_{w(2), n}\\ &&&I_{d_3}&\ast & \vdots\\ & & \Ddots && & \vdots\\ &I_{d_{r-1}} &\ast &\cdots&\ast & y_{w(n-1), n} \\ q &qy_{12} & & \cdots& &q y_{1n} \end{matrix}\right). $$
%The star-ed entries in the second matrix the
%Here all star-ed entries are some of the $x_{ij}$ or $y_{ij}$. 
Since $\gamma\in \Gamma_0(q)$, we must have $y_{12}, \ldots, y_{1, n-1} \in \Bbb{Z},$ hence by our choice of representatives
$$y_{12} =  \ldots =  y_{1, n-1}= 0.$$
Next we consider the ($n-1$)-st row 
$$(\underbrace{qx_{n, n-1}}_{\in \Bbb{Z}}, \underbrace{\ldots}_{d_{r-1} \text{ entries}}, \ast  %+ \underbrace{qy_{1, d+1} x_{n, n -1}}_{\in \Bbb{Z}}
, \cdots,  \ast %+ \underbrace{q y_{1, n-1} x_{n, n-1}}_{\in \Bbb{Z}}
, y_{w(n-1), n} + qy_{1n}x_{n-1, n}) \in \Bbb{Z}^n$$ of $\gamma$, where the stars are the same as the stars in the $(n-1)$-st row of $c^{\ast}wy$. We conclude that all star-ed entries in the $(n-1)$-st row of $c^{\ast} w y$   must be integral, hence 0.  %, i.e.\ $0$ modulo $\Bbb{Z}$. 
We continue with the $(n-2)$-nd row of $\gamma$. By the same argument we first have $qx_{n-2, n} \in \Bbb{Z}$ and then also star-ed entry in $(n-2)$-nd row of $x$ is integral (hence 0) as well as   all star-ed entries in the $(n-2)$-nd row of $c^{\ast} w y$.    Continuing in this way, all star-ed entries must vanish. In other words,  $x \, \text{diag}(1/q, 1\ldots, 1, q)\,  w y \in   \Gamma_0(q)$ with $x \in U(\Bbb{Q})$, $y \in U_w(\Bbb{Q})$ and all relevant entries in $[0, 1)$ implies
$$x = \left(\begin{matrix}  1 &  &\cdots    &x_{1}/q\\   & \ddots&  & \vdots\\    && 1& x_{n-1}/q \\   & & &1  \end{matrix}\right) , \quad y = \left(\begin{matrix}  1 &  &\cdots    &y_{1}/q\\   & \ddots&  & \vdots\\    && 1& y_{n-1}/q \\   & & &1  \end{matrix}\right) $$
with $x_i, y_i \in \{0, \ldots,  q-1\}$.  For these $x, y$ we compute
$$xc^{\ast} w y  = \left(\begin{matrix} x_{1} & & & &  (x_{1}y_{1} + 1)/q\\  x_{2 } & && I_2 &  (x_{2} y_{1} + y_{w(2)})/q\\ \vdots&& \Ddots&& \vdots \\  x_{n-1 }& I_{r-1}&&&  (x_{n-1 } y_{ 1} + y_{w(n-1) })/q\\ q &&&& y_{n}
\end{matrix}\right).$$
Obviously, this matrix is in $\Gamma_0(q)$ if and only if  the $n-1$ congruences
$$x_1y_1 + 1 \equiv 0 \, (\text{mod } q), \quad x_i y_1 + y_{w(i)} \equiv 0 \, (\text{mod } q), \quad 2 \leq i \leq n-1$$
are satisfied. This can be solved easily, and we obtain the explicit expression
$$S^v_{q, w}(M, N, (q, \ldots, q)) = \sum_{\substack{x_1, \ldots, x_{n-1} \, (\text{mod } q)\\ (x_1, q) = 1}} e\left(\frac{M_1 x_{n-1} \pm N_1 \bar{x}_1 x_{w^{-1}(n-1) } }{q}\right).$$
If $(M_1N_1, q) = 1$, the sum vanishes unless $n-1 = w^{-1}(n-1)$. The latter case happens for $w$ of the form \eqref{w}  with $d_1 = d_r = 1$, if and only if $w = w_{\ast}$ and then the Kloosterman sum equals $q^{n-2}$. \hfill $\square$

  \section{Automorphic forms and Whittaker functions}

 We denote by $\{\varpi\}$ an orthonormal basis  of right ${\rm O}_n(\Bbb{R}){\rm Z}^+$-invariant automorphic forms for the group $\Gamma_0(q)$, cuspidal or Eisenstein series.  The space $L^2(\Gamma_0(q)\backslash \mathcal{H})$ is equipped   with   the standard   inner product $\langle f, g \rangle = \int_{\Gamma_0(q)\backslash \mathcal{H}} f(xy)  \bar{g}(xy)  \dd x\, \dd^{\ast}y$. 
 We denote by $\int_{(q)} \dd \varpi$ a combined sum/integral over the complete spectrum of $L^2(\Gamma_0(q)\backslash \mathcal{H})$. The relevant spectral decomposition is a special case of Langlands' general theory, see e.g.\ \cite{Ar} for a convenient summary in adelic language. All $\varpi$ belong to representations of level $q' \mid q$ (cf.\ \cite[Th\'eor\`eme]{JPSS}) and we assume that  $\{\varpi\}$  contains all cuspidal  newvectors  of level $q' \mid q$. The underlying representation is denoted by $\pi$, so $\varpi \in V_{\pi}$. For notational simplicity let us denote the local archimedean Langlands parameter $\mu_{\pi}(\infty)$ simply by $\mu = (\mu_1, \ldots, \mu_n)$; it  satisfies
 %By Hecke theory, we can and will assume that all $\varpi$ are eigenfunctions of the (unramified) Hecke algebra coprime to $q$ and that $\{\varpi\}$ contains all cuspidal newforms. 
% At an unramified place $v$ the local Langlands parameter $\mu_{\pi}(v)$ satisfies
% We denote by $\mu = \mu_{\pi} \in \Bbb{C}^{n}$ with 
  \begin{equation}\label{unit}
     \mu_1 + \ldots + \mu_n = 0, \quad \{\mu_1, \ldots, \mu_n\} = \{-\bar{\mu}_1, \ldots ,- \bar{\mu}_n\}. 
     \end{equation}
 For a (not necessarily cuspidal) automorphic form $\varpi $ and   $N \in \Bbb{N}^{n-1}$ we define its $N$-th Fourier coefficient $A_{\varpi}(N)$ by
  \begin{equation}\label{four}
  \int_{U(\Bbb{Z})\backslash U(\Bbb{R})} \varpi(xy) \theta_N(-x) \dd x = \frac{A_{\varpi}(N)}{N^{\eta}} W_{\mu}(N\cdot  {\rm y}(y))
  \end{equation}
 where $y \in \tilde{T}(\Bbb{R})$ and $W_{\mu} : \Bbb{R}^{n-1}_{>0} \rightarrow \Bbb{C}$ is the standard (spherical) Whittaker function, cf.\ e.g.\ \cite[Section 2]{St}.
 
   If $\varpi$ is a cuspidal newform and $(m, q) = 1$, the $(m, 1,   \ldots, 1)$-th Fourier coefficient is proportional to the   $m$-th Hecke eigenvalue $\lambda_{\pi}(m)$ (or $\lambda_{\tilde{\pi}}(m)$ depending on normalization), and by  Rankin-Selberg theory  we obtain
\begin{equation}\label{Li}
|A_{\varpi}((m, 1, \ldots, 1))|^2 \asymp_{\mu} \frac{|\lambda_{\pi}(m)|^2 }{[{\rm SL}_n(\Bbb{Z}) : \Gamma_0(q)] L(1, \pi, \text{Ad}) } \gg_{\mu}  |\lambda_{\pi}(m) |^2 q^{-(n-1) -\varepsilon}
\end{equation}
if $\varpi$ is $L^2$-normalized, cf.\ e.g.\ \cite[Proposition 1]{Ve}.   Here we used the upper bound \cite[Theorem 2]{Li} for the $L$-value (the residue of the Rankin-Selberg $L$-function) on the edge of the critical strip. 

The following easy, but important lemma shows that $\lambda_{\pi}(p^{\nu})$ is (perhaps not always, but sufficiently often) as big as $p^{\nu \sigma_{\pi}(p)}$ with $\sigma_{\pi}(p)$ as in \eqref{sigma}. 

\begin{lemma}\label{Lenstra} For a prime $p \nmid q$ and $\nu > n$ we have
 $$\max_{0 \leq j \leq n-1} | \lambda_{\pi}(p^{\nu-j})| \geq (2 p^{\sigma_{\pi}(v)})^{1-n}   p^{\nu \sigma_{\pi}(p)}.$$
\end{lemma} 

\textbf{Proof.} The following argument is taken from   \cite[Lemma 3]{LS}. We have an identity of power series 
\begin{equation*}%\label{power}
\sum_{\nu = 0}^{\infty} \lambda_{\pi}(p^{\nu}) x^{\nu}   = \prod_{j=1}^n (1 - p^{\mu_{\pi}(p, j)} x)^{-1}.
\end{equation*}
%We have
%$$\lambda_{\pi}(p) = \sum_{j=1}^n p^{\mu_{\pi, j}(p)}.$$
Without loss of generality let $\mu_{\pi}(p, 1)$ have the largest real part, i.e.\ $\Re \mu_{\pi}(p, 1) = \sigma_{\pi}(p)$. Then
$$\sum_{\nu = 0}^{\infty} p^{\nu \mu_{\pi}(p, 1)} x^{\nu} = \prod_{j=2}^n (1 - p^{\mu_{\pi}(p, j)} x)  \sum_{\nu = 0}^{\infty} \lambda_{\pi}(p^{\nu}) x^{\nu} . $$
Comparing coefficients, we obtain the lemma. \hfill $\square$\\
%so that the sequence $(\lambda_{\pi}(p^{\nu}))_{\nu}$ satisfies a linear $n$-term recurrence with constant coefficients and at least one out of $n$ consecutive terms must grow according to a largest root:

We need an archimedean analogue of this result, which is a bit more technical. Roughly speaking, the growth of  $W_{\mu}$ near the origin should capture the size of $\sigma_{\pi}(\infty)$ in the same way as the growth of $\lambda_{\pi}(p^{\nu})$ captures the size of $\sigma_{\pi}(p)$, but this is harder to see as the Mellin transform of $W_{\mu}$ is not perfectly understood and the location of poles is subtle. We start by summarizing some properties.   As in \cite[(3.1), (3.2)]{St} we consider  the re-normalized Whittaker function
 \begin{equation}\label{normalize}
 W_{\mu}^{\ast}(y) = \pi^{(n-1)n(n+1)/12} y^{-\eta/2} W_{2\mu}\Big((\sqrt{y_1}/\pi, \ldots, \sqrt{y_{n-1}}/\pi)\Big).
 \end{equation}
 The corresponding Mellin transform $\widehat{W}^{\ast}_{\mu}(s) =  \int_{\Bbb{R}^{n-1}_{>0}} W^{\ast}_{\mu}(y)  y^{s } \frac{dy_1}{y_1} \cdots \frac{dy_{n-1}}{y_{n-1}} $ 
 is meromorphic in $\mu$ and $s \in \Bbb{C}^{n-1}$ \cite{FG}. Explicitly, we have \cite[(3.7)]{St}
 \begin{displaymath}
 \begin{split}
 &\widehat{W}^{\ast}_{\mu}(s_1) = \Gamma(s_1 + \mu_1)\Gamma(s_1+\mu_2), \quad n=2,\\
&  \widehat{W}^{\ast}_{\mu}(s_1, s_2) = \frac{\Gamma(s_1 + \mu_1)\Gamma(s_1+\mu_2)\Gamma(s_1 + \mu_3)\Gamma(s_2-\mu_1)\Gamma(s_2 - \mu_2)\Gamma(s_2+\mu_2)}{\Gamma(s_1+s_2)}, \quad n=3,
 \end{split}
 \end{displaymath}
 but in general there do not seem to be such simple formulae.   For $\Re s_2, \ldots, \Re s_{n-1}$ sufficiently large and $\Re s_1 > \sigma_{\pi}(\infty)$, the function $\widehat{W}^{\ast}_{\mu}(s)$ is holomorphic by \cite[Theorem 3.1]{St}. If in addition $\mu_1, \ldots, \mu_n$ are pairwise distinct, then $\widehat{W}^{\ast}_{\mu}(s)$ has  simple poles at $s_1 = -\mu_j$, $1 \leq j \leq n$, with residue
 $$\widehat{W}^{\ast}_{\mu^{(j)}}(s^{(j)})\prod_{\substack{1 \leq k \leq n\\ k \not= j}} \Gamma( \mu_k - \mu_j) $$
 where $$ \textstyle s^{(j)} =  (s_2, \ldots, s_{n-1}) + (\frac{n-2}{n-1}, \ldots, \frac{1}{n-1})\mu_j , \quad \mu^{(j)} = (\mu_1, \ldots, \mu_{j-1} , \mu_{j+1}, \ldots, \mu_{n-1}) +  \frac{\mu_j}{n-1} \cdot \textbf{1},$$ see \cite[Theorem 3.2]{St}. These statements are proved by a recursion formula \cite[(3.5)]{St}   of the form
 $$\widehat{W}^{\ast}_{\mu}(s) = \int \cdots \int \widehat{W}_{\nu}^{\ast}\Big(-t_1 - \frac{\alpha_1 + \alpha_2}{n-2}, \underbrace{\ast, \ldots, \ast}_{n-4}\Big) \Gamma(t_1 + s_1) (\ast) \frac{\dd t_1 \cdots \dd t_{n-3}}{(2\pi i)^{n-3}}$$
 where $\alpha_1, \alpha_2$ are any two elements from the multi-set $\{\mu_1, \ldots, \mu_{n}\}$ and  $\nu - \frac{\alpha_1 + \alpha_2}{n-2} \cdot \textbf{1} \in \Bbb{C}^{n-2}$ is the $(n-2)$-tuple of the remaining $\mu_j$; moreover,   $(\ast)$ is  independent of $s_1$ and holomorphic in $t_1$ in a wide vertical strip if $\Re s_2, \ldots, \Re s_{n-1}$ are sufficiently large, and the other $n-4$ arguments of $\widehat{W}_{\nu}^{\ast} $ are independent of $t_1$ and $s_1$. Inductively, starting from the explicit formula for $n=2$ and $n=3$, we see that in any fixed vertical strip for $s_1$ and for $\Re s_2, \ldots, \Re s_{n-1}$ sufficiently large, the only poles can occur at $s_1 = -\mu_j - k$ for $1 \leq j \leq n$, $k \in \Bbb{N}_0$. We conclude that
 $$\widehat{W}^{\dagger}_{\mu}(s) := \widehat{W}^{\ast}_{\mu}(s) \prod_{j=1}^n (s_1 + \mu_j)$$
is holomorphic for $\Re s_1 > \sigma_{\pi}(\infty) - 1$ (for sufficiently large $\Re s_2, \ldots, \Re s_{n-1}$) and
\begin{equation}\label{res}
\widehat{W}^{\dagger}_{\mu}(-\mu_j, s_2, \ldots, s_{n-1}) = \widehat{W}^{\ast}_{\mu^{(j)}}(s^{(j)})\prod_{\substack{1 \leq k \leq n\\ k \not= j}} \Gamma(1 + \mu_k - \mu_j) .
\end{equation}
 For this statement the assumption that the $\mu_j$ are pairwise distinct can be dropped by holomorphic continuation (note that by the Luo-Rudnick-Sarnak bounds or even the Jacquet-Shalika bounds $|\Re \mu_j| < 1/2$ the gamma factors on the right hand side  are always defined).

 For $\beta \in \Bbb{C}$ let $\mathcal{D}_{\beta} = -y\partial_{y} + \beta$. This is a commutative family of differential operators that under Mellin transformation correspond to multiplication with $s + \beta$. In the proof of Lemma \ref{whit} below we will need to following technical, but elementary lemma. 
 
 \begin{lemma}\label{new} Let $ \alpha \geq 0$, $c_0, c_1, c_2 > 0$,  $\beta \in \Bbb{C}$.  Let $I = [a, b]\subseteq  ( 0, 1)$ be an interval with $(1+c_0)a \leq b \leq 2a$ and $w : I \rightarrow \Bbb{C}$ a smooth function satisfying 
 \begin{equation}\label{diff}
    |\mathcal{D}_{\beta} w(y) | \geq c_1 y^{-\alpha}, \quad | \partial_y (\mathcal{D}_{\beta} w)(y) | \leq c_2 \| \mathcal{D}_{\beta} w \| y^{-1}
    \end{equation} for $y \in I$. Then there exist  constants $c_0',  c_1', c_2'> 0$ depending only on $c_0, c_1, c_2,   \beta$  (but not on $a, b$) and an interval $I' = [a', b']  \subseteq I$ with  $(b' - a') \geq c_0'(b-a)$ such that %such that $a \leq a' < b' \leq b$ and 
    \begin{equation}\label{diff1}
    |w(y)| \geq c_1' y^{-\alpha}, \quad |w'(y)| \leq c'_2 \| w_{I'}  \|y^{-1}
    \end{equation}
     for $y \in I'$. 
  \end{lemma}
  
  \textbf{Proof.} Let $\tilde{w}(y) = w(y) y^{-\beta}$, so that
  $$y^{1+\beta}\tilde{w}'(y) = -\mathcal{D}_{\beta}w(y), \quad y^{1+\beta}\tilde{w}''(y) = - \partial_y (\mathcal{D}_{\beta}w) (y) + \frac{1+\beta}{y} \mathcal{D}_{\beta}w(y).$$
  Then \eqref{diff} implies 
  $$|y \tilde{w}'(y)| \geq c_1 y^{-\tilde{\alpha}}, \quad  |\tilde{w}''(y)|  \leq \tilde{c}_2 \| \tilde{w}' \|y^{-1}$$
 for   $\tilde{c}_2 = 2^{1+|\Re \beta|} c_2 + |1+\beta|$ and $\tilde{\alpha} = \alpha + \Re \beta$.    Let $y_0 = \max_{y \in I} |\tilde{w}'(y)|$. Changing $\tilde{w}$ by a fourth root of unity if necessary, we can assume that $$\Re \tilde{w}'(y_0) \geq \frac{1}{\sqrt{2}} \max \big( c_1 y_0^{-\tilde{\alpha} - 1},   \| \tilde{w}' \| \big).$$
 The condition $|\tilde{w}''(y)| \leq \tilde{c}_2 \| \tilde{w}' \| y^{-1}$ implies that the slightly weaker inequality
 $$\Re \tilde{w}'(y) \geq \frac{1}{2\sqrt{2}} \max \big( c_1 y_0^{-\tilde{\alpha} - 1},   \| \tilde{w}' \| \big) \asymp \Re \tilde{w}'(y_0)$$
 holds on   some non-empty sub-interval $I_0 = [a_0, b_0] \subseteq I$ containing $y_0$, where $(b_0 - a_0) \geq (b-a)/\sqrt{8} \tilde{c}_2$.  Distinguishing the cases $\Re \tilde{w}(a_0) > - c_3 y_0 \Re \tilde{w}'(y_0)$ and $\Re \tilde{w}(a_0) \leq - c_3 y_0 \Re \tilde{w}'(y_0)$ for   $c_3 = (b_0-a_0)/\sqrt{8}y_0 \geq c_0 /32 \tilde{c}_2$, we confirm in both cases 
\begin{displaymath}
\begin{split}
%\| \tilde{w} \|_{I_0}&  \geq
 \max(|\Re\tilde{w}(a_0)|, |\Re\tilde{w}(b_0)|) \geq \min\Big(c_3 y_0 \Re \tilde{w}'(y_0), \Big(-c_3 y_0 + \frac{b_0-a_0}{2\sqrt{2}}\Big)\Re \tilde{w}'(y_0)\Big) = c_3  y_0 \Re \tilde{w}'(y_0).  %\geq \frac{c_3}{2\sqrt{2}} y |\Re \tilde{w}'(y)|
\end{split}
\end{displaymath}
 % $$y|\tilde{w}'(y)| \leq 2 \sqrt{2} y_0 \Re \tilde{w}'(y_0) \leq \min\Big(\frac{2\sqrt{2}}{c_3}, \frac{b_0 - a_0 - 2\sqrt{2}c_3}{}  \|\tilde{w} \| := \sup_{y\in I'} |\tilde{w}(y)|$$ on a suitable non-empty subinterval $I' \subseteq I_0$ 
 Thus there exists a non-empty subinterval $I' = [a', b'] \subseteq I_0$ of length $\gg (b_0 - a_0)$ such that 
% for $y \in I_0$ 
% which implies automatically 
$|\tilde{w}(y)| \geq \frac{1}{10} c_3 c_1 y^{-\tilde{\alpha}}$ on $I'$. Changing back to $w$, we obtain \eqref{diff1}. \hfill $\square$ \\
  
 We are now prepared for the following analogue of Lemma \ref{Lenstra}. For a function $E $ on $ \Bbb{R}^{n-1}_{>0}$ and $X\in \Bbb{R}_{>0}^{n-1}$ define 
 \begin{equation}\label{ex}
   E^{(X)}(y_1, \ldots, y_{n-1}) = E(X_1y_1,  \ldots, X_{n-1}y_{n-1}).
   \end{equation}
 
 \begin{lemma}\label{whit} Assume that $\mu$ varies in some compact set $\Omega$, and let $Z \geq 1$. 
% {\rm (a)} 
 There exist $r \in \Bbb{N}$ and a compact set $S\subseteq  \Bbb{R}_{>0}^{n-1}$ both depending only on $\Omega$ (not on $Z$) and a finite collection of (measurable) functions $E_1, \ldots, E_r : \Bbb{R}_{>0}^{n-1}\rightarrow \Bbb{R}$  h depending on $\Omega$ and $Z$ that are uniformly bounded (independent of $Z$) and supported in  
 $S$  such that 
  $$\sum_{j=1}^r |\langle E_j^{(Z, 1,\ldots, 1)}, W_{\mu}\rangle|^2 \gg_{\Omega} Z^{2\eta_1 +2 \sigma_{\pi}(\infty)}$$ for $\mu \in \Omega$ and $\eta$ as in \eqref{eta}. 
% {\rm (b)} Let $X > 1$ be a sufficiently large parameter. There exists a finite collection of  compactly supported functions $E_1, \ldots, E_r : \Bbb{R}_{>0}^{n-1}\rightarrow \Bbb{R}$, depending only on $\Omega$, such that $\sum_j |\langle E_j, W_{\mu}\rangle|^2 \gg 1$ for $\mu \in \Omega$.
  \end{lemma}
 
 \textbf{Proof.} For $Z \ll 1$ this is \cite[Lemma 1]{BBM}. For convenience we repeat the short argument in a slightly modified fashion that we will need later. There exists $Z_0 > 0$ (depending only on $\Omega$) such that for each $\mu \in \Omega$ we can choose an open set $S_{\mu} \subseteq \Bbb{R}_{>0}^{n-1}$ such that $|\Re W_{\mu}(y)|  > 2Z_0$ for all $y \in S_{\mu}$ or  $|\Im W_{\mu}(y)|  > 2Z_0$  for all $y \in S_{\mu}$.  Next  choose open neighbourhoods $U_{\mu}$ about $\mu$ such that  $|\Re W_{\mu^{\ast}}(y)|  > Z_0$ for all $y \in S_{\mu}$ and all $\mu^{\ast} \in U_{\mu}$ or  $|\Im W_{\mu}(y)| > Z_0$  for all $y \in S_{\mu}$ and all $\mu^{\ast} \in U_{\mu}$. By compactness we pick a finite collection  of such neighbourhoods $U_{\mu_1}, \ldots, U_{\mu_r}$ covering $\Omega$, and define the corresponding $E_j$ to be real-valued functions with support on $S_{\mu_j}$ and non-vanishing on the interior $\mathring{S}_{\mu_j}$. 
 
Now suppose that $Z$ is sufficiently large (in terms of $\Omega$). We try to mimic the proof of Lemma \ref{Lenstra}.  
%For $\beta \in \Bbb{C}$ let $\mathcal{D}_{\beta} = -y_1\partial_{y_1} + \beta$. This is a commutative family of differential operators that under Mellin transformation correspond to multiplication with $s_1 + \beta$. 
Assume (without loss of generality by \eqref{unit} and Weyl group symmetry) that $\Re (-\mu_1) = \sigma_{\pi}(\infty)$,  and with the notation as above let $$\widehat{{\tt W}}_{\mu}(s) := \frac{\widehat{W}^{\dagger}_{\mu}(s) }{s + \mu_1} = \widehat{W}^{\ast}_{\mu}(s) \prod_{j=2}^n (s_1 + \mu_j). $$
Taking inverse Mellin transforms, we obtain
$${\tt W}_{\mu}(y) = \mathcal{D}_{\mu_2} \cdots \mathcal{D}_{\mu_n} W_{\mu}^{\ast}(y)$$
where the differential operators are applied to the first variable $y_1$. 
 On the other hand, by Mellin inversion and \eqref{res} we have the asymptotic expansion
 $$y_1^j \partial_{y_1}^j{\tt W}_{\mu}(y) =  \mu_1^jy_1^{\mu_1}W^{\ast\ast}_{\mu}(y_2, \ldots, y_{n-1}) + O_{y_2, \ldots, y_{n-1}, \mu}\big(y_1^{\Re \mu_1 + 1/2}\big)$$
for $y_1 \rightarrow 0$ and $j \in \{0, 1\}$ where
$$W^{\ast\ast}_{\mu}(y_2, \ldots, y_{n-1}) =  W^{\ast}_{\mu^{(1)}}(y_2, \ldots, y_{n-1}) \prod_{j=2}^{n-1} y_j^{\frac{n-j}{n-1}\mu_1} \prod_{k=2}^n\Gamma(1 + \mu_k - \mu_1).$$
Whenever $|W^{\ast\ast}_{\mu}(y_2, \ldots, y_{n-1}) | \geq Z^{-1/2}$, say, and $y_2, \ldots, y_n \asymp 1$ (with implied constants depending only on $\Omega$), we can apply repeatedly   Lemma \ref{new} to ${\tt W}_{\mu}(y)$ with $\beta = \mu_2, \ldots, \mu_n$ to obtain two constants $1/2 < \gamma_1 < \gamma_2 < 1$ with 
%It is easy to see that if $\mathcal{D}_{\beta}w(y) \sim c y^{-a}$, $y \rightarrow 0$,  for some constants $a \geq 0$, $\beta, c \in \Bbb{C}$ and $Z\geq 1$ is sufficiently large, then there exist  constants $0 < \gamma_1 < \gamma_2 < 1$ (depending on all parameters, but uniformly bounded away from 0 as long as long as $\beta, c, a$ vary in a fixed compact set and $Z$ is sufficiently large) such that $|w(y)| \gg |c| y^{-a}$ for $y \in [\gamma_1/Z^2, \gamma_2/Z^2]$. Iterating this argument and adjusting the constants $\gamma_1, \gamma_2$ if necessary, we see that  
$$|W_{\mu}^{\ast}(y) | \gg  y_1^{-\sigma_{\pi}(\infty)} |W^{\ast\ast}_{\mu}(y_2, \ldots, y_{n-1})|$$ for   $y_1 \in [\gamma_1/Z^2, \gamma_2/Z^2]$. % as long as $\mu, y_2, \ldots, y_{n-1}$ vary in some fixed compact domain. 
By the same argument as in the beginning of the proof, we can now choose a finite collection of  functions $E_j^{\ast\ast} : \Bbb{R}_{>0}^{n-2} \rightarrow \Bbb{C}$ depending on $\Omega$ (but not on $Z$, provided that $Z$ is sufficiently large in terms of $\Omega$) such that $\sum_j |\langle E^{\ast\ast}_j, W^{\ast\ast}_{\mu} \rangle|^2 \gg 1$ for $\mu \in \Omega$, the inner product being restricted to the last $n-2$ coordinates.  Next define $E^{\ast}_j(y_1, \ldots, y_n) = \delta_{\gamma_1 \leq y_1 \leq \gamma_2} E^{\ast\ast}_j(y_2, \ldots, y_{n-1})$, %This choice depends on $Z$, but the support of $E$ varies inside some   interval depending only on $\Omega$. 
so that
  $$\sum_j \Bigl|\int_{\Bbb{R}_{>0}^{n-1}}   E^{\ast}_j(Z^2 y_1, y_2, \ldots, y_{n-1}) \overline{W^{\ast}_{\mu}(y)}   \frac{dy_1}{y_1} \cdots \frac{dy_{n-1}}{y_{n-1}} \Bigr|^2 \gg Z^{4\sigma_{\pi}(\infty)}.$$ Finally changing variables $y_j \leftarrow y_j^{1/2}/ \pi$ as in \eqref{normalize}, we obtain
\begin{displaymath}
\begin{split}
Z^{4\sigma_{\pi}(\infty)}&  \ll \sum_j  \Bigl|\int_{\Bbb{R}_{>0}^{n-1}} y^{-\eta} E_j^{\ast}(Z^2\pi^2 y_1^2, \pi^2 y_2^2, \ldots, \pi^2 y_{n-1}^2) \overline{W_{2\mu}(y)}   \frac{dy_1}{y_1} \cdots \frac{dy_{n-1}}{y_{n-1}} \Bigr|^2\\
& = Z^{-2\eta_1} \sum_j  \bigl|\langle E^{(Z, 1, \ldots, 1)}_j, W_{2\mu}\rangle \big|^2
\end{split}
\end{displaymath}
upon defining 
$E_j(y_1, \ldots, y_{n-1}) = y^{\eta} E_j^{\ast}(\pi^2 y_1^2, \pi^2 y_2^2, \ldots, \pi^2 y_{n-1}^2)$. Re-normalizing $\mu$ and $\sigma_{\pi}(\infty)$ by division by 2, we obtain the lemma. \hfill $\square$
%\int_{\Bbb{R}_{>0}^{n-1}} y^{-\eta} E^{\ast}(Z\pi^2 y_1^2, \pi^2 y_2^2, \ldots, \pi^2 y_{n-1}^2) \overline{W_{2\mu}(y)} y^{-2\eta} \frac{dy_1}{y_1} \cdots \frac{dy_{n-1}}{y_{n-1}} \Bigr|^2.$$

 \section{Poincar\'e series and the Kuznetsov formula}
  
   Let $E$ be a fixed compactly supported (measurable) function on $\Bbb{R}_{>0}^{n-1}$, $X \in \Bbb{R}_{>0}^{n-1}$ a ``parameter'' and define the right ${\rm O}_n(\Bbb{R}){\rm Z}^+$  invariant function $F^{(X)} : {\rm GL}_n(\Bbb{R}) \rightarrow \Bbb{C}$ by 
   \begin{equation}\label{F}
     F^{(X)}(xyk\alpha) = \theta(x) E^{(X)}({\rm y}(y)) %\quad E^{(X)}(z) = E(X\cdot z),
   \end{equation}
      for $x \in U(\Bbb{R})$, $y \in \tilde{T}(\Bbb{R})$, $k \in {\rm O}_n(\Bbb{R})$, $\alpha \in {\rm Z}^{+} $ and $\theta = \theta_{(1, \ldots, 1)}$ as in \eqref{character}, $E^{(X)}$ as in \eqref{ex}.     For $N \in \Bbb{N}^{n-1}$ we consider the Poincar\'e series
$$P^{(X)}_{N}(xy) = \sum_{\gamma \in U(\Bbb{Z}) \backslash \Gamma_0(q)} F^{(X)}(\iota(N)\gamma xy).$$
Note that $F^{(X)}(\iota(N)xy) = \theta_{N}(x) E( X\cdot N \cdot{\rm y}(y))$, cf.\ \eqref{char}.  % = \theta(\iota(N)x\iota(N)^{-1})E(\iota(N)y)$. 
Let $N, M \in \Bbb{N}^{n-1}$. By \cite[Theorem A]{Fr} (with $\rho = \text{triv}$, $\nu_1 = \ldots = \nu_{n-1} = 0$) we have
 \begin{equation}\label{four-poin}
 \begin{split}
& \int_{U(\Bbb{Z})\backslash U(\Bbb{R})} P^{(X)}_{M}(xy) \theta_N(-x) \dd x\\
 &  = \sum_{w\in W} \sum_{v \in V} \sum_{c \in \Bbb{N}^{n-1}} S_{q, w}^v( M, N, c) \int_{U_w(\Bbb{R}) } F^{(X)}(\iota(M)c^{\ast}w xy )\theta_N^v(-x) \dd x.
 \end{split}
 \end{equation}
 For fixed $y$ and fixed compact support of $E$, it follows from the two bounds in Lemma \ref{volume} with $M$ in place of $B$ that the $c$-sum runs over a finite set (depending on $M$, $y$ and the support of $F$), and  the $U_w(\Bbb{R})$-integral  runs over a compact domain (again depending on $M$, $y$ and the support of $F$). In particular the right hand side is absolutely convergent (and   the assumption $\Re \nu_j > 2/n$ in \cite[Theorem A]{Fr} can be dropped; Friedberg works more generally with bounded $E$ rather than compactly supported $E$). Without loss of generality we can assume that $w$ is of the form \eqref{w}.  

Now let $\varpi$ be a not necessarily cuspidal automorphic form   occurring in the spectrum of $L^2(\Gamma_0(q)\backslash \mathcal{H})$. By unfolding, \eqref{four} and a change of variables $y \leftarrow \iota(N)y$, we have
$$\langle  \varpi, P^{(X)}_{N}\rangle = \int_{\tilde{T}(\Bbb{R})} \int_{U(\Bbb{Z})\backslash U(\Bbb{R})} \varpi(xy) \theta_{ N}(-x) \overline{E^{(X)}(N \cdot {\rm y}( y))} \dd x\, \dd^{\ast} y  = N^{\eta}  A_{\varpi}( N) \langle W_{\mu}, E^{(X)}\rangle$$
where as before $\mu = \mu_{\pi}(\infty)$. 
By Parseval we obtain
$$\langle P^{(X)}_{M}, P^{(X)}_{N} \rangle = N^{\eta} M^{\eta} \int_{(q)} \overline{A_{\varpi}(M)}  A_{\varpi}(N) |\langle W_{\mu}, E^{(X)}\rangle|^2 \dd \varpi.$$
On the other hand,  by unfolding and \eqref{four-poin} we can express $\langle P^{(X)}_{M}, P^{(X)}_{N} \rangle$ as 
\begin{displaymath}
\begin{split}
%\langle P^{(X)}_{M}, P^{(X)}_{N} \rangle 
& \int_{\tilde{T}(\Bbb{R})} \int_{U(\Bbb{Z})\backslash U(\Bbb{R})} P^{(X)}_{M}(xy) \theta_{N}(-x) \overline{E^{(X)}(N \cdot {\rm y}( y))} \dd x\, \dd^{\ast}y\\
& =  \sum_{w\in W} \sum_{v \in V} \sum_{c \in \Bbb{N}^{n-1}} S_{q, w}^v( M, N, c)  \int_{\tilde{T}(\Bbb{R})} \int_{U_w(\Bbb{R}) } F^{(X)}(\iota(M)c^{\ast}w xy ) \theta_N^v(-x) \overline{E(X \cdot N \cdot {\rm y}( y))} \dd x\, \dd^{\ast}y.
\end{split}
\end{displaymath}
Let 
\begin{equation}\label{A}
A = \iota(X \cdot M)c^{\ast} w \iota(X \cdot N)^{-1} w^{-1} = \iota\big( X \cdot M \cdot{}^w(X \cdot N) \big)c^{\ast} \in T(\Bbb{R}), 
\end{equation}
so that ${\rm y}(A)^{\eta} c_1 \cdots c_{n-1} = \big( X \cdot M \cdot{}^w(X \cdot N) \big)^{\eta}$ by \eqref{ceta}. 
We change variables $y \leftarrow \iota(X \cdot N)y$, $x \leftarrow \iota(X \cdot N)x \iota(X \cdot N)^{-1}$. By Lemma \ref{lemma1}   we obtain
\begin{displaymath}
\begin{split}
&\sum_{w\in W} \sum_{v \in V} \sum_{c \in \Bbb{N}^{n-1}} S_{q, w}^v( M, N, c)  \frac{(X \cdot M)^{\eta} (X \cdot N)^{\eta}}{c_1 \cdots c_{n-1} {\rm y}(A)^{\eta}} \\
&\times \int_{\tilde{T}(\Bbb{R})} \int_{U_w(\Bbb{R}) } F^{(X)}(\iota(X)^{-1} A w xy ) \theta^v(-x) \overline{E( {\rm y}( y))} \dd x\, \dd^{\ast}y.
\end{split}
\end{displaymath}
%where
%By \eqref{yNw} we have 
%$${\rm y}\big(\iota(M)w \iota(N)^{-1} w^{-1} \big) = \Big(M_j\frac{N_1 \cdots N_{n - w(n-j+1)}}{N_1 \cdots N_{n- w(n-j)}}\Big)_{1 \leq j \leq n-1}.$$
We conclude the following Kuznetsov-type formula.

\begin{lemma}\label{kuz-formula} Let $M, N \in \Bbb{N}^{n-1}$, $X \in \Bbb{R}_{> 0}^{n-1}$, $E$ a compactly supported function on $\Bbb{R}_{> 0}^{n-1}$ and define $F^{(X)}$ as in \eqref{F}. Then
\begin{equation}\label{formula}
\begin{split}
&\int_{(q)} \overline{A_{\varpi}(M)}  A_{\varpi}(N)  |\langle W_{\mu}, E^{(X)}\rangle|^2 \dd \varpi \\
&= \sum_{w\in W} \sum_{v \in V} \sum_{c \in \Bbb{N}^{n-1}} \frac{S_{q, w}^v(M, N, c) }{c_1 \cdots c_{n-1}}  \frac{X^{2\eta}}{ {\rm y}(A)^{\eta}} \int_{\tilde{T}(\Bbb{R})} \int_{U_w(\Bbb{R}) } F^{(X)}( \iota(X)^{-1} A w xy ) \theta^v(-x) \overline{E( {\rm y}( y))} \dd x\, \dd^{\ast}y
\end{split}
\end{equation}
with $A$ as in \eqref{A}. 
\end{lemma}

As mentioned before, the Kloosterman sum $S_{q, w}^v(M, N, c)$ vanishes unless $w$ is of the form \eqref{w}, in which case we have the additional conditions \eqref{divi}, as well as \eqref{comp-equiv} for $i \not\in \{d_1, \ldots, d_1 + \ldots + d_{r-1}\}$. The $c$-sum is restricted by Lemma \ref{volume} and the support of $E$. %The volume of the $x$-integration is given in Lemma \ref{int-lemma}.  

\section{Proofs of Theorems \ref{thm1}, \ref{thm2}, \ref{thm4}}

We start with the \textbf{proof of Theorem \ref{thm2}}. We specialize Lemma \ref{kuz-formula} to $$M = N = (m,1,  \ldots, 1), \quad X = (Z, 1, \ldots, 1)$$ with $(m, q) = 1$.  We need to bound the spectral side from below and the Kloosterman side from above. By \eqref{Li} and positivity we have
\begin{equation}\label{11}
\sum_{\pi \in \mathcal{F}_I(q)} |\lambda_{\pi}(m)|^2  Z^{2\eta_1 + 2\sigma_{\pi}(\infty)} \ll_I q^{n-1+\varepsilon} \int_{(q)} |A_{\varpi}(M)|^2  Z^{2\eta_1 + 2\sigma_{\pi}(\infty)}   \delta_{\lambda_{\varpi} \in I} \,  \dd \varpi.
\end{equation}
By Lemma \ref{whit}  there is a finite set of compactly supported functions $E_j$ such that  $$Z^{2\eta_1 + 2\sigma_{\pi}(\infty)}  \delta_{\lambda_{\varpi} \in I} \ll_I \sum_j |\langle  W_{\mu_{\varpi}}, E^{(X)}_j\rangle|^2  .$$ 
Thus in order to bound the left hand side of \eqref{11} it suffices to  consider the right hand side of \eqref{formula} for a fixed $E^{(X)} = E^{(X)}_j$, and we are left with bounding
\begin{displaymath}
\begin{split}
q^{n-1+\varepsilon} \sum_{w\in W} \sum_{v \in V} &\sum_{c \in \Bbb{N}^{n-1}} \frac{S_{q, w}^v(M, N, c) }{c_1 \cdots c_{n-1}}  \frac{X^{2\eta}}{ {\rm y}(A)^{\eta}} \\
& \times \int_{\tilde{T}(\Bbb{R})} \int_{U_w(\Bbb{R}) } F^{(X)}( \iota(X)^{-1} A w xy ) \theta^v(-x) \overline{E( {\rm y}( y))} \dd x\, \dd^{\ast}y.
\end{split}
\end{displaymath}

For  $w = \text{id}$ we have  $c_1 = \ldots =c_{n-1} = 1$ and hence $A = I_n$, and the contribution is $  O(q^{n-1+\varepsilon}Z^{2\eta_1})$. 

Let us now consider the remaining $w$ of the form \eqref{w}. First we bound the  moduli $c_j$. To this end we apply   Lemma \ref{volume} with   $B = X \cdot M \cdot {}^w(X \cdot N)$, so that by \eqref{yNw} we obtain $B_1 = B_{n-d_1} = mZ$, $B_{n-d_1+1} = 1/(mZ)$ if $d_1 > 1$ and $B_j=1$ for all other indices. This gives 
\begin{equation}\label{cj}c_j \ll (mZ)^{s(1, j) + s(n-d_1, j) - s(n - d_1 + 1, j)} = \begin{cases} mZ, & j \leq n - d_1,\\
1, & j > n - d_1. \end{cases}
\end{equation}
 We assume that $mZ \ll q^2$ with a sufficiently small implied constant, so that $c_j < q^2$ for all $j$. We may also assume that $q$ is sufficiently large, otherwise there is nothing to prove. Now suppose that $d_1 > 1$ (but $d_1 < n$ since $w \not= \text{id}$). Then by \eqref{comp-equiv} with $i =  d_1 - 1$ we have $c_{n-d_1 + 2} c_{n- d_1} = \pm c_{n- d_1 + 1}^2 m$. Using \eqref{divi} and comparing the $q$-adic valuation on both sides, we conclude from \eqref{divi} that both $c_{n-d_1 + 1}$ and $c_{n- d_1 + 2}$ are divisible by $q$, which contradicts  \eqref{cj} for $q$ sufficiently large. Hence $d_1 = 1$, and we see from \eqref{divi} that all $c_j$ are divisible by $q$. We write $c_j = qc_j'$. By \eqref{mult}
 % Applying the same argument for the condition \eqref{comp-equiv} with $i =   d_1 - 2$, we see that also $q\mid c_{n-d_1 + 3}$, and inductively all $c_j$ are divisible by $q$, exacly once. But then \eqref{comp-equiv} with $i =1$ gives $c_{n-2} = \pm c_{n-1}^2$, contradiction.  We conclude that $d_1 = 1$. Similarly and even more easily, if $d_r > 1$, we have 
we obtain
\begin{equation}\label{chin}
S_{q, w}^v(M, M, c) = S^v_{q, w}(\ast, \ast, (q, \ldots, q)) S_{1, w}^v(M, (\bar{q}m, 1, \ldots, 1, \bar{q}), c')
\end{equation}
where   $\ast$ is coprime to $q$. By Theorem \ref{thm3} the first factor on the right hand side vanishes unless  $w = w_{\ast}$, in other words, only the trivial Weyl element and $w_{\ast}$ survive.  

(As an aside: if we only wanted to prove Sarnak's density original density hypothesis with an exponent $n- 1 - 2\sigma + \varepsilon$ in Theorem \ref{thm1}, then upon choosing $mZ \ll q$ with a sufficiently small constant, \emph{all} Weyl elements except the identity would vanish and no further analysis would be necessary. That in the stronger set-up $mZ \ll q^2$ only $w_{\ast}$ needs to be considered is an artefact of $q$ being prime.)

Our next aim is to estimate 
\begin{displaymath}
\begin{split}
&\Big| \int_{\tilde{T}(\Bbb{R})} \int_{U_{w_{\ast}}(\Bbb{R}) } F^{(X)}(\iota(X)^{-1}A w_{\ast} xy ) \theta^v(-x) \overline{E( {\rm y}( y))} \dd x\, \dd^{\ast}y\Big| \\
&\leq  \int_{\tilde{T}(\Bbb{R})} \int_{U_{w_{\ast}}(\Bbb{R}) } |E({\rm y}(A w_{\ast} xy) )  E( {\rm y}( y))| \dd x\, \dd^{\ast}y.
\end{split}
\end{displaymath}
 By  Lemma  \ref{volume} and then    Lemma \ref{int-lemma}   the right hand side is bounded by 
\begin{equation}\label{72}
\begin{split}
 %&  \int_{\tilde{T}(\Bbb{R})} \int_{U_{w_{\ast}}(\Bbb{R}) } F^{(X)}(\iota(X)^{-1}A w_{\ast} xy ) \theta^v(-x) \overline{E( {\rm y}( y))} \dd x\, \dd^{\ast}y \\
 &\ll_E \text{{\rm vol}}\Big\{x \in U_{w_\ast}(\Bbb{R}) \mid \Delta_j(w_{\ast}x) \ll_E \prod_{i=1}^{n-1} {\rm y}(A)_i^{s(i, j)}, 1 \leq j \leq n-1\Big\} \\
 &\ll_E \prod_{i=1}^{n-1}\prod_{j=1}^{n-1} {\rm y}(A)_i^{s(i, j)(1+\varepsilon)} = {\rm y}(A)^{\eta(1+\varepsilon)}
   \end{split}
\end{equation}
since $\sum_{i} s(i, j) = \eta_j$ by \eqref{inverse} and \eqref{eta}. 

Summarizing the previous  estimations (and changing the value of $\varepsilon$) and applying Theorem \ref{thm3} to the first factor in \eqref{chin}, we obtain
\begin{displaymath}
\begin{split}
&\sum_{\pi \in \mathcal{F}_I(q)} |\lambda_{\pi}(m)|^2  Z^{2\eta_1 + 2\sigma_{\pi}(\infty)} \\
&\ll_I Z^{2\eta_1} q^{n-1+\varepsilon} \Big(1 + q^{n-2} \sum_{v \in V} \sum_{c'_1, \ldots, c_{n-1}' \ll mZ/q} \frac{|S^v_{1, w_{\ast}}(M, (\bar{q}m, 1, \ldots, 1, \bar{q}),  c')|}{q^{n-1} c_1' \cdots c_{n-1}'}\Big).
\end{split}
\end{displaymath}
For the Weyl element $w_{\ast}$ the consistency relations \eqref{comp-equiv} impose serious restrictions on the moduli $c_1', \ldots, c_{n-1}'$. We apply \eqref{comp-equiv} with $i = 2, \ldots, n-2$ getting 
$(c'_i)^2 = c'_{i-1}c'_{i+1}$ for $i = 2, \ldots, n-2$. %, which implies
%$$c'_i (c'_{n-1})^{n-i-2} = (c'_{n-2})^{n-i-1}, \quad i = 1, \ldots, n-3,$$
%as one easily checks.
If $n \geq 4$, then $c'_{2}$ fixes $c_1'$ and $c_3'$ up to a divisor function, and inductively also $c_4', \ldots, c_{n-1}'$.   Using the trivial bound \eqref{triv}, we finally obtain (again changing the value of $\varepsilon$)
$$\sum_{\pi \in \mathcal{F}_I(q)} |\lambda_{\pi}(m)|^2 Z^{ 2\sigma_{\pi}(\infty)}  \ll_I q^{n-1+\varepsilon} \Big(1 + \frac{q^{n-2}}{q^{n-1}}   \sum_{ c_{2}' \ll mZ/q} 1  \Big) \ll q^{n-1+\varepsilon}$$
provided $mZ \ll q^2$. In the case $n=3$ we quote from \cite[(4.2) with $N=1$]{BBM} the average Weil-type bound
$$\sum_{c_1, c_2 \leq X} |S^v_{w^{\ast}}((m, 1), (\bar{q}m,  \bar{q}),  c')|  \ll X^3 (Xm)^{\varepsilon}$$
to obtain again
$$\sum_{\pi \in \mathcal{F}_I(q)} |\lambda_{\pi}(m)|^2 Z^{ 2\sigma_{\pi}(\infty)}  \ll_I q^{2+\varepsilon} \Big(1 + \frac{1}{q} \cdot \frac{mZ}{q}  \Big) \ll q^{2+\varepsilon}$$
for $mZ \ll q^2$. This completes the proof. \hfill $\square$ \\

The \textbf{proof of Theorem \ref{thm4}} is a simple variation. %Here we do not assume that $q$ is prime. 
Again by   positivity and \eqref{Li} we have
\begin{displaymath}
\begin{split}
 \sum_{\pi \in \mathcal{F}_I(q)} & \Big| \sum_{\substack{m \leq x\\ (m, q) = 1}} \alpha(m) \lambda_{\pi}(m) \Big|^2 \ll q^{n-1+\varepsilon}   \int_{(q)}  \Big| \sum_{\substack{m \leq x\\ (m, q) = 1}} \alpha(m)A_{\varpi}(M) \Big|^2     \delta_{\lambda_{\varpi} \in I} \,  \dd \varpi\\
 & = q^{n-1+\varepsilon}  \sum_{\substack{m_1, m_2 \leq x\\ (m_1m_2, q) = 1}}\alpha(m_1) \overline{\alpha(m_2)}   \int_{(q)}    A_{\varpi}(M_1) \overline{A_{\varpi}(M_2)}      \delta_{\lambda_{\varpi} \in I} \,  \dd \varpi
 \end{split}
\end{displaymath}
where $M = (m, 1, \ldots, 1)$, $M_1 = (m_1, 1, \ldots, 1)$, $M_2 = (m_2, 1, \ldots, 1)$. 
We detect the condition $\delta_{\lambda_{\varpi} \in I}$ by a finite collection of test functions $E_j$ with $Z = 1$ as in the previous proof and apply Lemma \ref{kuz-formula}. For $w \not= \text{id}$ the analogue of \eqref{cj} is 
$$c_j \ll m_2^{s(1, j) } m_1^{s(n-d_1, j) - s(n - d_1 + 1, j)} \leq x$$
which contradicts \eqref{divi} for $x \ll q$ (with a sufficiently small implied constant) since $d_1 \not= n$. So only the trivial Weyl element survives, and we obtain the desired bound. \hfill $\square$\\

\textbf{Corollary \ref{cor5}} follows easily Theorem \ref{thm4} by observing that an approximate functional equation has length $q^{1/2}$ (see \cite[Section 5]{IK}): for all but $O(1)$ cuspidal representations $\pi \in \mathcal{F}_I(q)$ (and $\varepsilon < 1/2$)  we have 
$$|L(1/2 + it, \pi)|^2 \ll_{I, t, n, \varepsilon} q^{\varepsilon}  \sum_{ 2^j = M \leq q^{1/2 + \varepsilon}  } \frac{1}{M}  \Bigl| \sum_{ M \leq m \leq 2M }  \lambda_{\pi}(m) \Bigr|^2 $$ %\ll q^{n-1 + \varepsilon}$$
and the desired bound follows directly from Theorem \ref{thm4}. Note that    the shape of the ramified coefficients (i.e.\ $q \mid m$) is irrelevant and the condition $(m, q) = 1$ in Theorem \ref{thm4} is void.  \\ 
 
Finally we derive \textbf{Theorem \ref{thm1}} from Theorem \ref{thm2}. Let us first assume that $v = p\not= q$ is a fixed prime.   %, as an implied constant on the right hand side one can choose $(2 p^{\sigma_{\pi}(v)})^{1-n}$. 
We choose $\nu_0$ maximal so that $ p^{\nu_0} \ll q^2$ with an  implied constant that is admissible for Theorem \ref{thm2}. %By Lemma \ref{Lenstra} we have 
%Next we choose $\nu_0 - n \leq \nu \leq \nu_0$ such that 
%$$\sum_{\nu_0 - n \leq \nu \leq \nu_0} |\lambda_{\pi}(p^{\nu})|^2 \gg p^{\nu \sigma_{\pi}(p)}.$$% according to Lemma \ref{Lenstra}.  
%Then clearly $p^{\nu} \asymp q^2$, and 
We conclude from Theorem \ref{thm2} with $Z = 1$, $m = p^{\nu}$  and Lemma \ref{Lenstra} that
\begin{displaymath}
N_p(\sigma, \mathcal{F}_I(q)) \leq \sum_{\pi \in \mathcal{F}_I(q)} \frac{p^{2\nu_0 \sigma_{\pi}(p)}}{p^{2\nu_0 \sigma}} \ll  \frac{1}{q^{4\sigma}} \sum_{\nu_0 - n \leq \nu \leq \nu_0} \sum_{\pi \in \mathcal{F}_I(q)}  |\lambda_{\pi}(p^{\nu})|^2  \ll q^{n-1-4\sigma + \varepsilon}.
  \end{displaymath}
For $v = \infty$, Theorem \ref{thm1} follows directly from Theorem \ref{thm2} with $m=1$, $Z \ll q^2$ (again with a sufficiently small implied constant).  \hfill $\square$.

\end{document}